\documentclass[a4paper]{article}

\usepackage[margin=1in]{geometry} 

\usepackage{amsthm}
\usepackage{amssymb}
\usepackage{amsmath}
  \usepackage{paralist}
  \usepackage{graphics} 
  \usepackage{epsfig} 
\usepackage{graphicx,color}
\graphicspath{{Figure/}}
\usepackage{epstopdf}
\usepackage{dsfont}
\usepackage{caption}
\usepackage{subcaption}
\usepackage{array,multirow,makecell}
\usepackage{booktabs}
\usepackage{lmodern}
\usepackage{arydshln}
\usepackage[utf8]{inputenc}
\usepackage{hyperref}

\usepackage{biblatex} 
\theoremstyle{plain}
\newtheorem{theorem}{Theorem}

\newtheorem{lemma}[theorem]{Lemma}

\newtheorem{proposition}[theorem]{Proposition}

\theoremstyle{definition}

\newtheorem{remark}[theorem]{Remark}

\usepackage{mathrsfs} 

\newcommand{\dint}{\displaystyle{\int}}
\newcommand{\e}{\mbox{e}}
\newcommand{\ud}{\mbox{d}}
\newcommand{\dsum}{\displaystyle{\sum}}

\title{Analysis of a tumor growth model with a nonlocal boundary condition}
\author{Slah Eddin Ben Abdeljalil \and Atef Ben Essid  \and Saloua Mani Aouadi\thanks{Corresponding author: Saloua Mani Aouadi (saloua.mani@fst.utm.tn) }}

\date{
	Faculty of Sciences of Tunis, University Tunis EL Manar, 2092, Tunisia 
	}

\begin{document}
	\maketitle

		\begin{abstract}
In this paper, we conduct a thorough mathematical analysis of a tumor growth model with treatments. The model is a system describing the evolution of metastatic tumors and the number of cells present in a primary tumor. The former evolution is described by a linear transport equation and the latter by an ordinary differential equation of Gompertzian type. The two dynamics are coupled through a nonlocal boundary condition that takes into account the tumor colonization rate. We prove an existence result where the main difficulty is to deal with the coupling and to take into account the time discontinuities generated by treatment terms. We also present numerical tests that highlight the effect of different treatments.

	\end{abstract}

	
\section{Introduction}\label{sec1}
In deciding the best treatment for cancer therapy, estimations of the colony size of tumors, predictions of the metastasis propagation and responses to treatments are needed. Given its vitality, this subject is of increasing interest for mathematicians \cite{literature,Luigi,Gandolfi,Cabrales,ghaffari}.
In this context, we consider the model of tumor growth introduced in \cite{IKS}. It is a linear coupling system between a transport partial differential equation (PDE) and an ordinary differential equation (ODE) with a non local boundary condition. We point out that, in \cite{Barbolosi}, the authors studied the existence, uniqueness, and asymptotic behavior of solution for this model by a semi-group approach. In \cite{FV,BS}, the authors added a term associated with chemotherapy, proved the existence of a solution, and carried out a thorough numerical analysis justifying the model by clinical trials. We note that in \cite{saut2017}, the authors proposed a PDE model describing cells movements generated by the proliferation of cancer cells which exert pressure on surrounding tissues. A PDE-ODE system modeling the multiscale invasion of tumor cells through the surrounding tissue matrix is studied in \cite{Christina2016,Christina2014}. In \cite{ghaffari}, some ODE models with chemo and radiotherapy treatments are analyzed. In \cite{KHAJANCHI,LI1,LI2}, the authors considered some conceptual models for the tumor-immune interaction based on dynamic systems, they focused on bifurcation and stability. The parameters' estimation, with the objective of obtaining prognostic models and of finding fields not observable, has been widely studied in the literature \cite{saut2012, Saut2014}. 

In this work, we complete the model analyzed in \cite{Barbolosi} by considering a concomitant treatment. The radiotherapy term is the one given in \cite{Radio1,Radio2}. The associated mathematical data are then discontinuous essentially in time. We focus on these irregularities in time and prove the existence and uniqueness of solution with piecewise regularity in time and Lebesgue integrability in space. Our strategy is the following: we first solve the ODE (decoupled part of the system) within the Caratheodory framework \cite{Daniel,Hauray}, then we plug the solution into the PDE part solving it using a fixed point argument. We also carry out some significant numerical tests. We discretize the ODE using Runge Kutta schemes and the PDE by the method of characteristics \cite{Bengt}. We notice that there is a large disparity coming from biologically relevant parameters, which is restrictive for the choice of the discretization parameters. Our numerical tests suggest that a mixed treatment should be better to reduce the population of cancer cells for some parameters related to the progress of the disease.
Finally, the contribution in this work lies at the same time in the modeling, the theoretical analysis, and the numerical approach. To the best of our knowledge, Iwata's model has not been used with mixed treatments. Also, compared to \cite{FV}, our proof of solution's existence is strongly penalized by the discontinuities of the data. On the numerical level, a sophisticated calculation code, based on the method of characteristics, is implemented. The numerical results show the interest in a mixed treatment for a clinical combination of parameters linked to the disease.

The paper is sketched as follows. In section 2, we write our model with some details, that allow us to understand the origin of irregularities. Section 3 is devoted to the analysis of the ODE's solution. In section 4, we present an existence result. The proof is based on the construction of an adequate contracting operator which takes into account the different constraints of the problem. Finally in section 5, we present our discretization strategy and some numerical tests.
\section{The mathematical model}
\subsection{Growth model and distribution of metastatic tumors}
\subsubsection{The primary tumor}
To formulate the process of metastases, we consider the model introduced in \cite{IKS} and taken up in \cite{Barbolosi,BS,ATP,hartung}. The authors consider an idealized case in which a primary tumor is generated from a single cell at time $t=0$ and grows at rate $g(x)$ per unit time, where $x$ denotes the tumor size represented by the number of cells in the tumor. The number $x_{p}(t)$ of cells present in a primary tumor at time $t$ is presented by the solution of
\begin{eqnarray}
\label{dyna1}
\frac{d}{dt}x_{p}(t) &=&g\left(x_{p}(t)\right);\text{ \ \ }t\geq 0, \\
 x_{p}(0) &=&1,  \notag
\end{eqnarray} 
with a Gompertzian growth rate
\begin{equation}
	g(x)=ax\ln \left(\frac{b}{x}\right)
	\end{equation}%
where $b\geq 1$ is the maximal tumor size and $a$ is a positive parameter that measures the rate of clonogenic proliferation. The solution of (\ref{dyna1}) is 
\begin{equation}
x_{p}(t)=b^{1-\e^{-at}} .  \label{xp}
\end{equation}%
 It increases strictly between $1$ and $b$, the corresponding curve is a sigmoid with three phases: a slow growth phase, an exponential growth phase and a slowdown one. This takes into account the stages of carcinogenesis, namely: the initiation phase corresponding to the transformation of a normal cell into a malignant cell, the phase of promotion corresponding to the accumulation of mutations and the slowing phase which corresponds to cell loss.
\subsubsection{The model without treatment}
Following \cite{IKS}, the primary growing tumor emits metastatic cells at rate $\beta (x)$. In turn, each metastatic cell grows into a new tumor which grows at rate $g(x)$ and emits new metastatic nuclei like the primary tumor. Let $u(t,x)$ represents the colony size distribution of metastatic tumors with $x$ cells at time $t$. Assuming the localized colonization nuclei sufficiently distant from each other, so that they do not overlap, the dynamic model proposed in \cite{IKS} writes

\begin{eqnarray}
\frac{\partial u}{\partial t}(t,x)+\frac{\partial }{\partial x}\bigl[g(x)u(t,x)\bigr]
&=&0,\text{ \ \ \ \ \ }t>0,\text{ }x\in ]1,b[,  \label{MacVon1} \\
u(0,x) &=&0,  \label{MacVon2} \\
g(1)u(t,1) &=&\int_1^b\beta (x)u(t,x)\ud x+\beta \left(x_{p}(t)\right).
\label{MacVon3}
\end{eqnarray}
The term $\frac{\partial }{\partial x}(g(x)u(t,x))$ reflects the transport of tumor cells. Equation (\ref{MacVon2}) indicates that initially no metastatic tumor exists. Equation  (\ref{MacVon3}) means that the number of metastatic single cells newly created per unit time is the total rate of new metastases due to metastases already present and to the primary tumor. The colonization rate $\beta (x)$ has the form
\begin{equation}
\beta (x)=m.x^{\alpha },  \label{taux}
\end{equation}
where $m$ is the colonization coefficient and $\alpha $ is the fractal dimension of the blood vessels infiltrated into the tumor. Equation (\ref{taux}) indicates that the rate of metastases from a tumor of size $x$ is proportional to the number of tumor cells in contact with the blood vessels. The parameter $\alpha $ expresses how the blood vessels are geometrically distributed in the tumor. We refer to \cite{IKS} for more details on the development of the model and its justification by clinical trials, as well as for an extensive bibliography on the subject.
\subsubsection{Consideration of chemo and radiotherapy}
In \cite{FV}, the authors have taken up the model of Iwata et al \cite{IKS} by considering chemotherapy treatment. We complete here by taking into account concomitant chemo and radiotherapy treatments. The Gompertzian growth rate $g$ is replaced by
\begin{equation}
G(t,x)=g(x)-K_{c}(x)C(t)-K_{r}(x)R(t),  \label{G}
\end{equation}
where the term $K_{c}(x)C(t)$ is associated with chemotherapy, $C$ represents the drug's concentration and $K_{c}$ measures its efficiency. We adopt for $K_{c}$ the expression
introduced in \cite{FV} 
\begin{equation*}
K_{c}(x)=\gamma (x-\overline{x})H(x-\overline{x}),
\end{equation*}%
where $ H $ is the Heaviside function, $ \overline {x} $ is a threshold from which drugs start and $ \gamma $  is a positive constant that quantifies the drug's effectiveness. Likewise, $ R $ represents cell death induced by radiation, and
\begin{equation*}
   K_{r}(x)=\gamma _{r}(x-\hat{x})H(x-\hat{x}), 
\end{equation*}
indicates that radiotherapy is applied to tumors larger than $\hat{x}.$ The terms $C$, $R,$ and the concomitant timeline depend on the types of protocols.
\paragraph{Chemotherapy}
The term $C(t)$ depends on the pharmacokinetics of administered drugs. The toxic effects of treatment that induce the lag term in the equations here are not considered. These toxic effects are taken into account \cite{Meille2008,Faivre}. Here, we use a model with a central compartment and a peripheral one. Let $V$ be the volume distribution, $k_{e}$ the elimination constant and $d_{c}$ the rate of drug's infusion into the central compartment. The drugs concentration is modeled by \cite{Faivre,FV,BS}

\begin{equation}
\left\{ 
\begin{array}{lclr}
\dot{c_{1}}(t) & = & -k_{e}c_{1}(t)+k_{12}(c_{2}(t)-c_{1}(t))+\dfrac{d_{c}(t)%
}{V} & \forall t\geq t_{0},\,c_{1}(t_{0})=0, \\ 
\dot{c_{2}}(t) & = & k_{21}(c_{1}(t)-c_{2}(t)) & \forall t\geq
t_{0},\,c_{2}(t_{0})=0,%
\end{array}%
\right.  \label{C1}
\end{equation}%
where $ t_ {0} $ is the start time of treatment, $ c_ {1} $, $ c_ {2} $ represent respectively the evolution of drug's concentration in the central and peripheral compartments, $ k_ {12 } $, $ k_ {21} $
are exchange constants between the two compartments and 
\begin{equation}
d_{c}=\dsum_{i=1}^{n}\frac{d_{i-1}^{c}}{t_{i}-t_{i-1}}\chi _{\lbrack
t_{i-1},t_{i}]}, \label{dose}
\end{equation}%

 $d_{i}^{c}$ is the dose given during the time $[t_{i-1},t_{i}]$. For the following, we notice that the second member of the differential system (\ref{C1}) is discontinuous.
\paragraph{Radiotherapy}
Using a linear-quadratic model, the cell survival probability writes
\begin{equation}
\text{survival probability}=\exp \left(-\alpha _{\text{eff}}d_{r}\right),  \label{prob}
\end{equation}%
where $d_{r}$\ is a radiation dose, and $\alpha _{\text{eff}}$ is a constant that translates radiations into cell death, said relative effective radiosensitivity parameter
  \cite{Radio1,Radio2,livreRadio,Radio4,Radio5,Radio6,Radio7,Julien}. In general, fractions of doses with the same magnitude are administered. We denote by $D_{r}(t)$ the accumulated dose at time $t.$ By not considering any delayed or otherwise toxic effects, the probability of cells death by radiation at time $t$ writes 
\begin{equation*}
R(t)=1-\exp\left(-\alpha _{\text{eff}}D_{r}(t)\right).
\end{equation*}%
Each radiation session lasts a few minutes, a finite series of pulses is
then administered. We assume that during these minutes, the irradiation takes place continuously and uniformly, so that the accumulated dose writes
\begin{equation*}
D_{r}(t)=\sum_{i=1}^{m}\frac{d_{r}}{2\varepsilon }\chi _{\lbrack
t_{i-\varepsilon },t_{i+\varepsilon }]}(t),
\end{equation*}%
where $m,\varepsilon $ and $t_{i}$ are characteristics of each protocol. Therefore, $R$ writes 
\begin{equation*}
R(t)=1-\exp \left(-\alpha _{\text{eff}}\frac{d_{r}}{2\varepsilon }%
\sum_{i=1}^{m}\chi _{\lbrack t_{i-\varepsilon },t_{i+\varepsilon }]}(t)\right).
\end{equation*}

\subsection{Final model}
The model that we are studying, finally writes 
\begin{eqnarray}
\frac{\partial }{\partial t}u(t,x)+\frac{\partial }{\partial x}\bigl[G(t,x)u(t,x)\bigr] &=&0,\text{ \ \ \ \ \ }t\in ]t_{0},T[,\text{ \ }x\in ]1,b[,
\label{Eq1} \\
(Gu)(t,1) &=&\int_{1}^{b}\beta (x)u(t,x)\ud x+f(t),\text{ \ \ \ }t\in ]t_{0},T[,
\label{Eq2} \\
(Gu)(t,b) &=&0,\text{ \ \ \ \ \ }t\in ]t_{0},T[,  \label{Eq3} \\
u(t_{0},x) &=&u_{0}(x),  \label{Eq4} \\
f(t) &=&\beta \left(x_{p}(t)\right),  \label{Eq5} \\
\frac{dx_{p}}{dt}(t) &=&G\left(t,x_{p}(t)\right),\text{ \ \ }t>t_{0},\text{ }%
x_{p}(t_{0})=x_{0},  \label{Eq6}
\end{eqnarray}%
where $t_{0}$ and $T$ are the start and end times of treatment, $u_{0}$ is
the density in size at $t_{0}$ and $G$ is given by (\ref{G}).

The existence of a solution, for the above problem, was proven in \cite{FV} with
$G\in \ C^{2}([t_{0},T]\times \lbrack 1,b])$. But taking the treatments into account makes $ G $ discontinuous, we focus on the irregularities in time : $K_{c}(t,.)$ and $K_{r}(t,.)$ will be approached with functions of class $C^{2}([1,b])$ for all $t$. 
 
\textbf{Hypothesis H.}\\ 
We will assume, in a generic framework, that $G$ is of class $C^{2}$  in $x$ for all $t$, has a finite number of discontinuities denoted by $t_{i},\; i=1,...,n$, and is of class $C^{2}$  on the connected components of $\mathcal{Q}=[t_0,T]\times[1,b]$ delimited by the curves $\Sigma_{i}=\{(t_{i},x), x\in[1,b]\},\;  i=1,...,n$. These hypothesis will be called \textbf{Hypothesis H}.

\section{The ordinary differential equation}
In this section, we are interested in the ODE (\ref{Eq6}). With \textbf{Hypothesis H}, this ODE has no classical solution (in the Cauchy sense), we look then for solutions in the Caratheodory sense \cite{Daniel,Hauray}.
\begin{proposition}
\label{Caract} We suppose that $G$ verifies \textbf{Hypothesis H}. For all $%
(s,x)\in \mathcal{Q}=[t_{0},T]\times \lbrack 1,b],$ the differential equation 
\begin{equation}
\left\{ 
\begin{array}{c}
\frac{d}{dt}\Phi (t)=G\bigl(t,\Phi (t)\bigr)\text{ a.e. in time,} \\ 
\Phi (s)=x,%
\end{array}%
\right.  \label{flot}
\end{equation}%
has a unique solution of class $C^{1}$ piecewise on $[t_{0},T]$ given by 
\begin{equation*}
\Phi _{(s,x)}(t)=x+\int_{s}^{t}G\left(r,\Phi _{(s,x)}(r)\right)\ud r.  \label{flot2}
\end{equation*}%

Moreover, $\Phi :(t,s,x)\rightarrow \Phi _{(s,x)}(t)$  is continuous on $[t_{0},T]\times \lbrack t_{0},T]\times \left[ 1,b\right]$,\newline
$\tilde{\Phi}_{t}:(s,x)\rightarrow \Phi _{(s,x)}(t)$ is of class $C^{1}$ in a neighborhood of any point where $G$ is $C^{1}$ for all fixed $t$,\newline
and $\hat{\Phi}_{(t,s)}:x\rightarrow \Phi _{(s,x)}(t)$ is of class $ C^{1}$ on $]1,b[$ for any fixed $(t,s)$.
\end{proposition}
\begin{proof}
With \textbf{Hypothesis H}, the function $G$ verifies the following the Caratheodory's conditions \cite{Daniel}
\begin{equation}
\begin{array}{l}
1.\text{\ for almost all }t \in [t_{0},T],\text{ }G(t,.)\text{ is continuous in }x\text{,} \\
2.\text{ for almost all }x \in [1,b],\text{ }G(x,.)\text{ is Lebesgue measurable in }t,   \\
3.\text{ there is a function }m \in L^{1}([t_{0},T])\text{ such that} \\
\qquad\qquad\left\vert G(t,x)\right\vert \leq m(t)\text{ \ \ }\forall (t,x)\in Q, \\
4.\text{ there is a function }l \in L^{1}([t_{0},T])\text{ such that}\\
\qquad\qquad\left\vert G(t,x)-G(t,y)\right\vert \leq l(t)\left\vert x-y\right\vert\text{ \ \ }\forall t\in \lbrack t_{0},T],\text{ \ }\forall x,y\in \lbrack1,b],  
\end{array}
\end{equation}%
with $m(t)=m=\max_{(t,x)\in \mathcal{Q}}G(t,x)$ and $l(t)=\max_{x\in \lbrack 1,b]}\frac{\partial G(t,x)\text{\ }}{\partial x}.$ Then, the system (\ref{flot}) admits a unique global solution in $W^{1,1}([t_{0},T])$ given by the integral equation (\ref{flot2}). In
addition, $\frac{d}{dt}\Phi _{(s,x)}(.)$ is continuous for all $t\neq t_{i},i=1,...,n.$ Otherwise, for all $(t,s_{1},x_{1}),(t,s_{2},x_{2})$ in $[t_{0},T]^2\times \left[ 1,b\right]$, we have for $s_1<s_2<t$
\begin{eqnarray*}
\left\vert \Phi _{(s_{1},x_{1})}(t)-\Phi _{(s_{2},x_{2})}(t)\right\vert &\leq &\left\vert x_{1}-x_{2}\right\vert+\Big\vert \int_{s_{1}}^{t}G\left(\tau,\Phi _{(s_{1},x_{1})}(\tau )\right)d\tau \\&&-\int_{s_{2}}^{t}G\left(\tau ,\Phi_{(s_{2},x_{2})}(\tau )\right)\ud\tau \Big\vert \\
&\leq &\left\vert x_{1}-x_{2}\right\vert +\int_{s_{1}}^{s_{2}}\left\vert G\left(\tau ,\Phi _{(s_{1},x_{1})}(\tau )\right)\right\vert \ud\tau\\&&+\int_{s_{2}}^{t}\left\vert G\left(\tau ,\Phi _{(s_{2},x_{2})}(\tau )\right)-G\left(\tau,\Phi _{(s_{1},x_{1})}(\tau )\right)\right\vert \ud\tau \\
&\leq &\left\vert x_{1}-x_{2}\right\vert +m\left\vert s_{1}-s_{2}\right\vert\\&&+\dint_{s_{2}}^{t}\left\vert G\left(\tau ,\Phi _{(s_{2},x_{2})}(\tau )\right)-G\left(\tau,\Phi _{(s_{1},x_{1})}(\tau )\right)\right\vert \ud\tau \\
&\leq &\left\vert x_{1}-x_{2}\right\vert +m\left\vert s_{1}-s_{2}\right\vert\\&&+\left\Vert l\right\Vert _{L^{\infty }(t_{0},T)}\dint_{s_{2}}^{t}\left\vert
\Phi _{(s_{2},x_{2})}(\tau )-\Phi _{(s_{1},x_{1})}(\tau )\right\vert \ud\tau.
\end{eqnarray*}
Gronwall's Lemma implies that
\begin{eqnarray} 
\left\vert \Phi_{(s_{1},x_{1})}(t)-\Phi _{(s_{2},x_{2})}(t)\right\vert &\leq &(\left\vert x_{1}-x_{2}\right\vert +m\left\vert
s_{1}-s_{2}\right\vert )\e^{ \left\Vert l\right\Vert _{L^{\infty}(t_{0},T)}\left\vert t-s_{2}\right\vert }.
\label{eq1}
\end{eqnarray}%
Now, for all $(t_{1},s_{1},x_{1})$, $(t_{2},s_{2},x_{2})$ in $%
[t_{0},T]\times \lbrack t_{0},T]\times \left[ 1,b\right] ,$ we have
\begin{eqnarray}
\left\vert \Phi (t_{1},s_{1},x_{1})-\Phi (t_{2},s_{2},x_{2})\right\vert&\leq &\left\vert \Phi (t_{1},s_{1},x_{1})-\Phi
(t_{1},s_{2},x_{2})\right\vert\\&&+\left\vert \Phi (t_{1},s_{2},x_{2})-\Phi(t_{2},s_{2},x_{2})\right\vert,  \label{eq2} 
\end{eqnarray}%
hence, the function : $(t,s,x)\rightarrow \Phi (t,s,x)$ is continuous. The
regularity of $\tilde{\Phi}$\ and $\hat{\Phi}_{(t,s)}$\ derives from
Cauchy's theorem.
\end{proof}

We will prove, in the following theorem, that the solutions of (\ref{Eq6}) allow ${Q}=[t_{0},T]\times \lbrack 1,b]$ to be partitioned into three
areas, each area will be concerned only by one boundary condition in problem (\ref{Eq1})--(\ref{Eq6}). In particular, we construct regular functions corresponding to the entry of the characteristic curves in each zone.
\begin{theorem}
\label{pointsflot} We suppose that $G$ verifies \textbf{%
Hypothesis H} and 
\begin{eqnarray}
G(t,1) &>&0, \text{ }\forall t\in [t_{0},T],  \label{E1} \\
G(t_{0},x) &>&0, \text{ }\forall x\in \lbrack 1,b[,  \label{E2} \\
G(t,b) &<&0, \text{ }\forall t\in ] t_{0},T].  \label{E3}
\end{eqnarray}%
We introduce the sets
\begin{eqnarray*}
\Omega _{1} &=&\big\lbrace(t,x)\in \mathcal{Q}\,\,\text{\textit{/} }x<\Phi
_{(t_{0},1)}(t)\big\rbrace, \\
\Omega _{2} &=&\big\lbrace(t,x)\in \mathcal{Q}\,\,\text{\textit{/}}\,\,\Phi
_{(t_{0},1)}(t)<x<\Phi _{(t_{0},b)}(t)\big\rbrace, \\
\Omega _{3} &=&\big\lbrace(t,x)\in \mathcal{Q}\,\,\text{\textit{/ }}\,x>\Phi
_{(t_{0},b)}(t)\big\rbrace.
\end{eqnarray*}%
Then there are a piecewise $C^{1}$ function $\varphi : %
\bar{\Omega}_{1}\rightarrow \lbrack t_{0},T]$ defined by 
\begin{equation}
\Phi \bigl(t,\varphi (t,x),1\bigr)=x,  \label{FI}
\end{equation}%
a piecewise $C^{1}$ function $\psi :\bar{\Omega}%
_{2}\rightarrow ]1,b[$ defined by 
\begin{equation*}
\Phi \bigl(t,t_{0},\psi (t,x)\bigr)=x,
\end{equation*}%
and a piecewise $C^{1}$ function $\theta :\bar{\Omega}%
_{3}\rightarrow \lbrack t_{0},T]$ defined by 
\begin{equation*}
\Phi \bigl(t,\theta (t,x),b\bigr)=x.
\end{equation*}
\end{theorem}
\begin{remark}
The hypotheses (\ref{E1}), (\ref{E2}) and (\ref{E3}) are verified by the
mapping $G$ given by (\ref{G}). They make sure that $\Omega _{1},$ $%
\Omega _{2}$ and $\Omega _{3}$ cover $[t_{0},T]\times \lbrack
1,b]$ and that each subdomain is concerned by only one boundary condition in the problem (\ref{Eq1})--(\ref{Eq6}). 
\end{remark}
\begin{figure}[ht]
\centering
\includegraphics[scale=0.3]{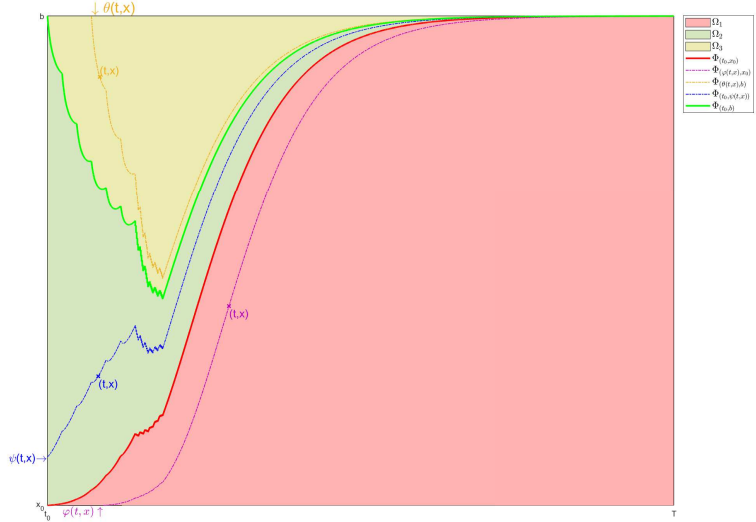} \vspace{-5mm}
\caption{The characteristic curves $\Phi $ on $\Omega _{1},\,\Omega _{2}$
and $\Omega _{3}$. }
\label{figure1}
\end{figure}
\begin{remark}
The functions $\varphi $ and $\theta $ define the entry time into $\Omega
_{1}$ (resp. $\Omega _{3}$) of the characteristic curves intersecting $\Omega _{1}$ (resp. $\Omega _{3}$). The function $\psi $ defines
the entry point in $\Omega _{2}$ of the characteristic curves intersecting $\Omega _{2}$ (see Figure. \ref{figure1}).
\end{remark}
\begin{proof}
\begin{enumerate}[1.]
\item We look for a time $\varphi (t,x)$ (resp. $\theta (t,x)$) such  that
\begin{eqnarray*}
\Phi \bigl(t,\varphi (t,x),1\bigr) &=&x\qquad (%
\text{resp.}\;\Phi \bigl(t,\theta (t,x),b\bigr)=x) ,\\
&&\text{or equivalently} \\
\Phi \bigl(\varphi (t,x),t,x\bigr) &=&1\qquad (\text{resp.}\;\Phi \bigl(\theta (t,x),t,x\bigr)=b).
\end{eqnarray*}
Let $(t,x)$ in $\Omega _{1}$. The characteristic curve resulting from $(t,x)$ intersects necessarily the line $x=1$ at a unique time $(\varphi (t,x),1)$ so that $(\varphi (t,x),1)$ is the foot of this characteristic curve. Likewise, for $(t,x)$ in $\Omega _{3}$, the  
characteristic curve passing through $(t,x)$ intersects necessarily the line 
$x=b$ at a unique time $\theta (t,x)$. This ensures the existence of continuous functions $\varphi $ and $\theta $.

Let us prove the regularity of $\varphi$ and $\theta$. 

For $(t,x)$ be in $\Omega _{1}$, $\varphi (t,x)$ is the entry time in $\Omega _{1}$. The mapping $G$ may be continuous or not in time at $\varphi(t,x)$.
\begin{enumerate}[a.]
\item If $G(.,x)$ is continuous at $s=\varphi (t,x)$, we have
\begin{equation*}
\Phi _{(t,x)}(s)=1\Rightarrow G\left(s,\Phi _{(t,x)}(s)\right)>0, \\ 
\end{equation*}%
and so $$\frac{\partial \Phi }{\partial s}(s,t,x)\neq 0.$$%
The function: $s\rightarrow \Phi (s,t,x)$ is then of class $C^{1}$ in a
neighborhood of $\varphi (t,x)$ and $\dfrac{\partial
\Phi }{\partial s}\bigl(\varphi (t,x),t,x\bigr)$ is strictly positive. We
deduce from the implicit function theorem that $s=\varphi (t,x)$ and that $\varphi $ is of class $C^{1}$ in a neighborhood of $(t,x)$.
\item If $G$ is not continuous in time at the point $\varphi (t,x)$ , the regularity of $\Phi $ implies that the curves $\Phi (s,t,x)$ and $\Phi (s,t^{\prime },x)$ are close and do not intersect
themselves for $t^{\prime }$ close to $t$ 
\begin{eqnarray*}
t\neq t^{\prime},\,\left\vert  t-t^{\prime}\right\vert<\varepsilon& \Rightarrow &\varphi (t,x) \neq  \varphi(t^{\prime },x)  \text{ and } \left\vert \varphi(t,x)-\varphi(t^{\prime},x)\right\vert\ll\varepsilon,\\
\end{eqnarray*}%
$\varphi (t^{\prime },x)$ is thus
a continuity point of $G(.,x)$, and $\varphi $ admits a right
and a left continuous derivative in time at the point $(t,x)$.

In the same way, with $x^{\prime }$ close to $x$ and the curves $\Phi
(s,t,x) $, $\Phi (s,t,x^{\prime }),$ we get that $\varphi $ is of class $C^{1}$ piecewise on $\Omega _{1}$.  

The exactly same approach with 
\begin{equation*}
\theta _{(t,x)}(s)=b\Rightarrow G\left(s,\Phi _{(t,x)}(s)\right)<0, \\ 
\end{equation*}%
ensures the regularity of $\theta$.
\end{enumerate}
\item For $(t,x)$ given in $]t_{0},T[\times ]1,b[$, we look for an
entry point $y\in ]1,b[$ that satisfies $\Phi (t,t_{0},y)=x.$ Immediately
we have $y=\Phi (t_{0},t,x)=\psi (t,x)$ and $\psi $ is of class $C^{1}$
piecewise in $\Omega _{2}$.
\end{enumerate}
\end{proof}
\begin{remark}
The mapping $\varphi $ is of class $C^{1}$ on parts of $\Omega _{1}$
delimited by the curves $\widetilde{\Sigma}_{i}=\Phi
_{(t_{i},1)}$, $i=1,...,n,$. The mapping $\psi $ is of class $C^{1}$ on the
parts of $\Omega _{2}$ delimited by $\Sigma _{i}=\{(t_{i},x),x\in
\lbrack 1,b]\},$ $i=1,...,n.$ And by construction, $\varphi $ and $\psi $
are constant along the characteristic curves so that 
\begin{equation*}
\varphi \bigl(s,\Phi (s,t,x)\bigr)=\varphi (t,x)\ ;\ \ \psi \bigl(s,\Phi (s,t,x)\bigr)=\psi
(t,x)\ \ \forall s,\text{ }t,\text{ }x.
\end{equation*}
\end{remark}
\begin{figure}[h!]
\begin{minipage}[c]{.49\linewidth}
\centering
\includegraphics[width=0.95\linewidth,height=0.8\linewidth]{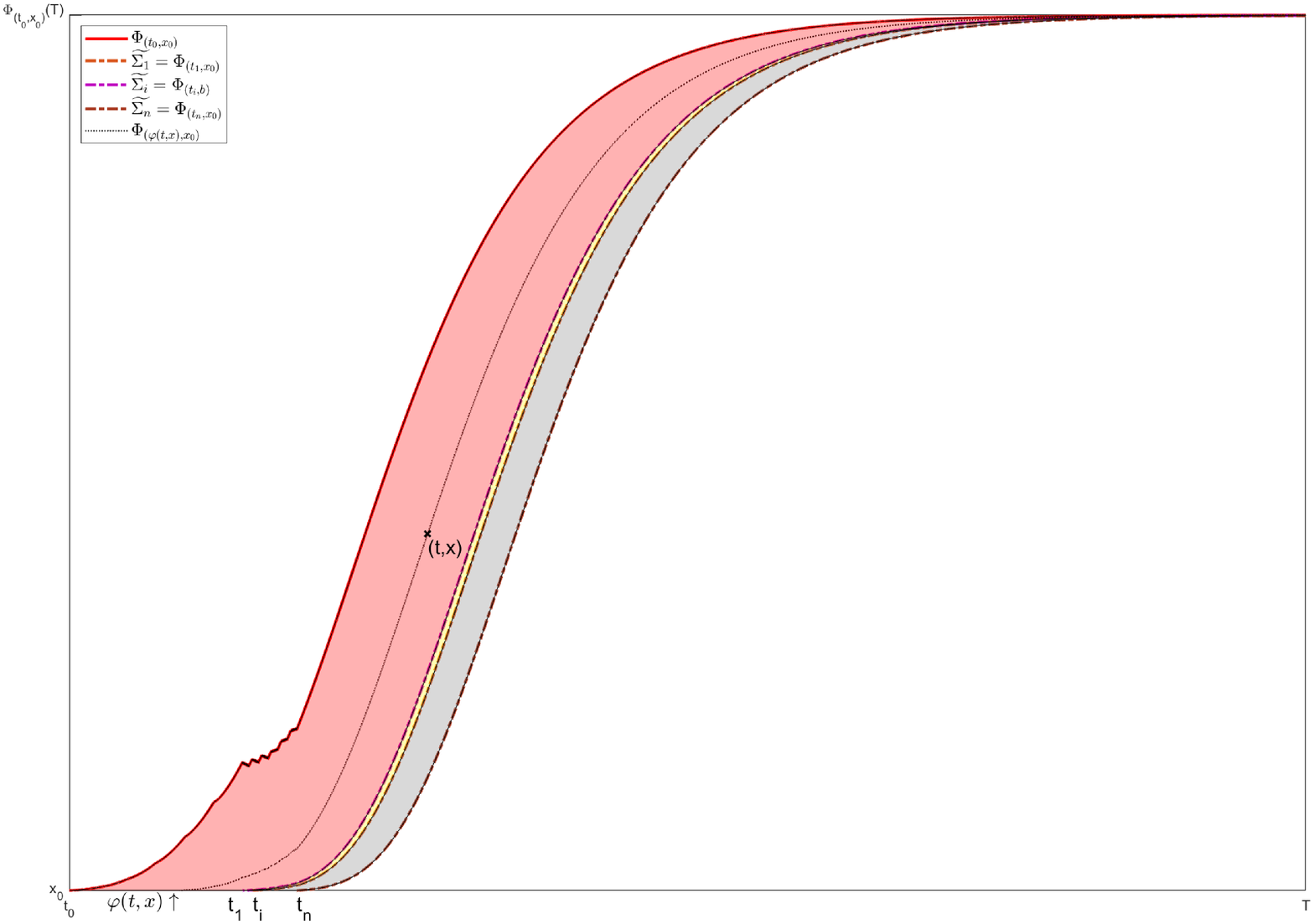}
\vspace{-3mm}
\caption[l]{Regularity domain of $\varphi$ delimited by $\widetilde{\Sigma_i}$.}
\label{figure2}
\end{minipage}
\hfill
\begin{minipage}[c]{.49\linewidth}
\centering
\includegraphics[width=0.95\linewidth,height=0.8\linewidth]{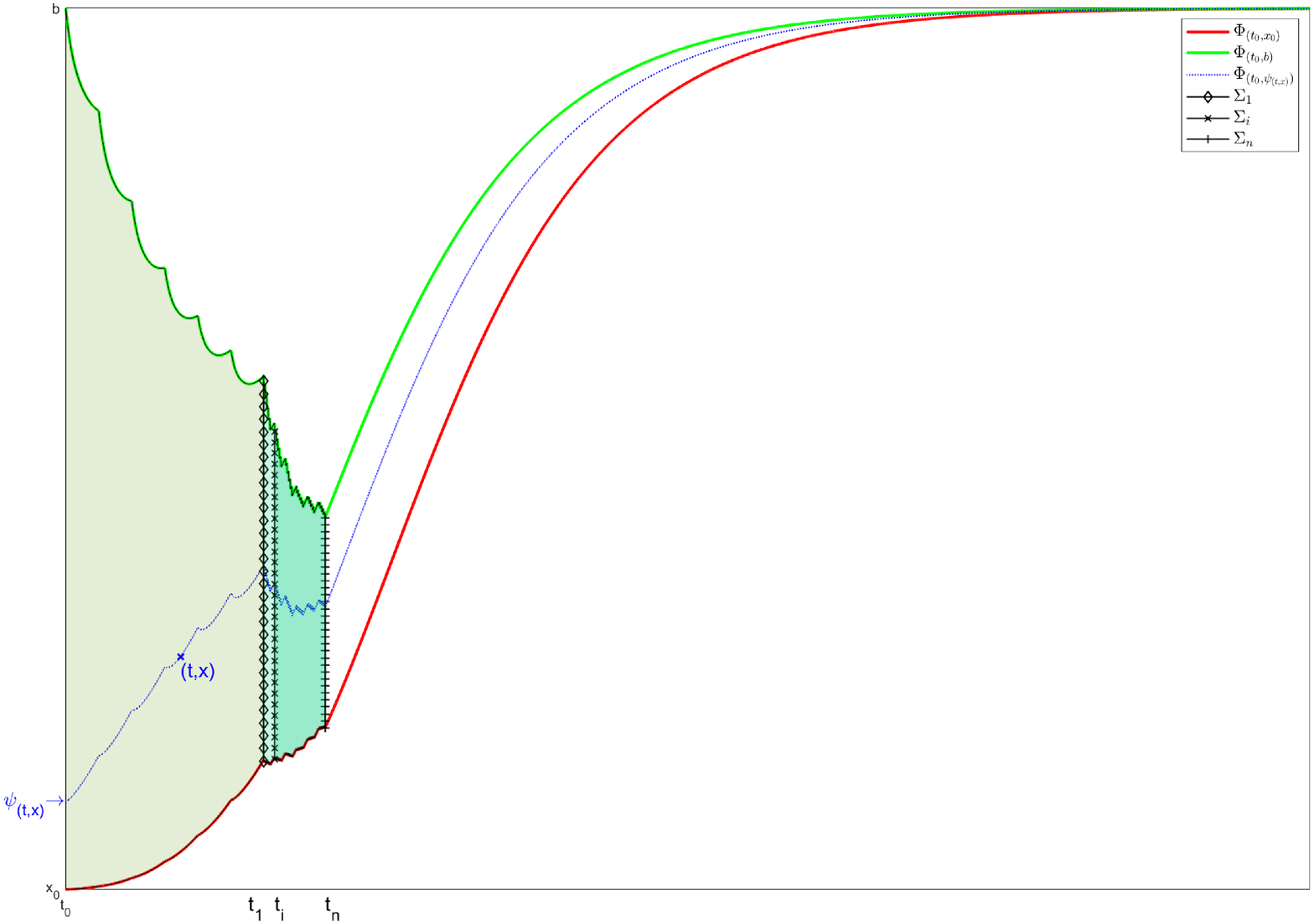}
\caption[l]{Regularity domain of $\psi$ delimited by $\Phi_{(t_0,b)}$, $\Phi_{(t_0,x_0)}$ and $\Sigma_i$.}
\label{figure3}
\end{minipage}
\end{figure}
\section{Existence result}
In this section, we focus on the proof of the existence and uniqueness of solution for
the problem (\ref{Eq1})-(\ref{Eq6}).
\subsection{Functional framework}
We denote by $E=\mathcal{C}_{MG}^{1}([t_{0},T],\,L^{1}([1,b]))$ the set of
continuous functions on $[t_{0},T]$ with values in $L^{1}([1,b])$, of class $%
C^{1}$ on $[t_{0},T]\backslash \{t_{i},\;i=1,...,n\}$ and admitting right
and left derivatives in time on each point. We take note that $E$ is different
from the set of class $C^{1}$ piecewise in time functions. The index $G$
reminds that the set of derivative discontinuity points is fixed, finite,
and is associated with the set of time discontinuities of $G$. We provide $E$
with its natural norm 
\begin{equation*}
\left\Vert w\right\Vert =\sup_{t\in \lbrack t_{0},T]}\left\Vert
w(t,.)\right\Vert _{L^{1}([1,b])}+\max_{i\in \{0,...,n\}}\sup_{t\in \lbrack
t_{i},t_{i+1}]}\left\Vert \frac{\partial w}{\partial t}(t,.)\right\Vert
_{L^{1}([1,b])}.
\end{equation*}%
We define, for $u\in L^{1}([1,b])$, the set 
\begin{equation*}
E(u)=\big\lbrace w\in \mathcal{C}_{MG}^{1}([t_{0},T],\,L^{1}([1,b]))\quad
w(t_{0},\,.)=u\big\rbrace.
\end{equation*}%
It is obvious that $(E,\left\Vert .\right\Vert )$ is a Banach space and that the metric
space $E(u)$ is closed in $E$ for each $u\in L^{1}([1,b])$.

\subsection{The main theorem}
\begin{theorem}
\label{th_existence_unicite}
We suppose that the initial data $u_{0}$ belongs to $W^{1,1}(]1,b[)$ and
verifies the compatibility conditions 
\begin{equation}
u_{0}(b)=0\quad\text{and}\quad(Gu_{0})(t_{0},1)=\int_{1}^{b}\beta
(x)u_{0}(x)\ud x+f(t_{0}).
\end{equation}
We suppose that $G$ verifies \textbf{Hypothesis H} and the conditions (\ref{E1})--(\ref{E3}) with $G(t,1)$ constant on $[t_{0},T]$. 

Then, there is a unique $u\in E(u_{0})$ with $Gu\in \mathcal{C}%
([t_{0},T],\,W^{1,1}([1,b]))$\ that verifies the equations (\ref{Eq1}), (\ref{Eq2}) and (\ref{Eq3}) respectively in $\mathcal{C}_{MG}%
([t_{0},T],\,L^{1}([1,b]))$, $\mathcal{C}([t_{0},T],\,\mathbb{R})$ and $%
L^{1}([1,b]).$
\end{theorem}
\begin{remark}
The hypothesis $G(t,1)$ constant is verified by $G$ given by (\ref{G}). It reflects the fact that radiotherapy starts only from a threshold $\hat{x}$ greater than $1$.
\end{remark}
The key ingredient in the proof of theorem \ref{th_existence_unicite} is a Banach fixed point theorem in a suitable fundamental functional space. The proof is deduced from the following technical four lemma. In the first one, we construct an operator on $E(u)$. In the second lemma, we prove that this operator has a fixed point. In the third lemma, we prove that this fixed point is a solution. In the last lemma, we prove the uniqueness.

\begin{lemma}
\label{const} Let $u_{0} \in W^{1,1}(]1,b[)$ given with $u_{0}(b)=0$. We
assume that the hypotheses of Theorem \ref{th_existence_unicite} are satisfied. We define the operator $\mathcal{T}$ by 
\begin{equation}
\mathcal{T}(w)(t,x) =u(t,x) =\left\{ 
\begin{array}{ll}
R\bigl(\varphi (t,x)\bigr)\exp \left(-\dint_{{\varphi (t,x)}}^{t}\partial _{x}G\left(s,\Phi
_{(t,x)}(s)\right)\,\mathrm{d}s\right) & \text{in }\Omega _{1}, \\ 
u_{0}\bigl(\psi (t,x)\bigr)\exp \left(-\dint_{t_{0}}^{t}\partial _{x}G\left(s,\Phi _{(t,x)}(s)\right)\,%
\mathrm{d}s\right) & \text{in } \Omega _{2}, \\ 
0 & \text{in } \Omega _{3},%
\end{array}%
\right. 
\label{Cons}
\end{equation}%
where 
\begin{equation*}
R(t)=\frac{1}{G(t,1)}\left( \int_{1}^{b}\beta (y)w(t,y)\,\mbox{d}%
y+f(t)\right).
\end{equation*}%
Then $\mathcal{T}$ is well defined from $E(u_{0})$ into $E(u_{0})$.
\end{lemma}
\begin{proof} Let $w$ be given in $E(u_{0})$. 
\begin{enumerate}[Step 1.]
\item The function $u$ defined by (\ref{Cons}) satisfies
\begin{eqnarray}
\dint_{1}^{b}|u(t,x)|\mathrm{d}x &\leq &M_{0}\bigg(C_{1}(\Vert \beta
\Vert _{\infty }\sup_{[t_{0},T]}\Vert w(t,.)\Vert _{L^{1}}  \notag +\sup_{t\in \lbrack t_{0},T]}|f(t)|)+\Vert u_{0}\Vert _{L^{1}}\bigg)<\infty
\label{majoration1}
\end{eqnarray}%
where 
\begin{equation}
M_{0}=\sup_{\mathcal{Q}}\exp \left( -\int_{t_{0}}^{t}\partial _{x}G\left(s,\Phi
_{(t,x)}(s)\right)\,\mathrm{d}s\right) \text{ and }C_{1}=\dfrac{b-1}{G(t,1)}.
\label{M0}
\end{equation}%
So, $u(t,.)$ is in $L^{1}(]1,b[).$ In addition, it is clear that $u$ is
continuous on ${Q}$ deprived of the characteristics curves $\Phi
_{(t_{0},1)} $ and $\Phi _{(t_{0},b)}.$ We also have, for any sequence $(\hat{t}%
_{n})_{n\in \mathbb{N}}$ that converges to $t$ 
\begin{eqnarray*}
\Vert u(\hat{t}_{n},.)-u(t,\,.)\Vert _{L^{1}} &\leq &M_{0}\Bigg[%
C_{1}\Vert \beta \Vert _{\infty }\Vert w(\varphi (\hat{t}%
_{n}),.)-w(\varphi (t),\,.)\Vert _{L^{1}}\\&& +C_{1}|f(\hat{t}_{n})-f(t)|+
\Vert u_{0}(\psi (\hat{t}_{n},.))-u_{0}(\psi
(t,\,.))\Vert _{L^{1}}\Bigg].
\end{eqnarray*}%
As $w\in \mathcal{C}([t_{0},T],$ $L^{1}(]1,b[)),$ $u_{0}\in \
W^{1,1}(]1,b[),\,f\in \mathcal{C}([t_{0},T])$ and $\varphi ,\,\psi $ are
continuous, we get 
\begin{equation*}
\Vert u(\hat{t}_{n},\,.)-u(t,\,.)\Vert _{L^{1}}\xrightarrow[{n\to+\infty}]{}0.
\end{equation*}
\item Let us prove that $u$ is continuous on the
characteristic curves $\Phi _{(t_{0},1)}$ and $\Phi _{(t_{0},b)}$. Let 
$(\tau ,x_{0})$ be a point on the curve $\Phi _{(t_{0},1)}.$ Due to
the continuity of $u_{0},$ $\varphi $ and the fact that ${\varphi (\tau ,x}%
_{0}{)=t}_{0}$, we have
\begin{equation*}
\begin{array}{ll}
\underset{(t,x)\in \Omega _{1}}{\underset{(t,x)\rightarrow (\tau ,x_{0})}{%
\lim }}u(t,x) & =R\bigl(\varphi (\tau ,x_{0}{)\bigr)\exp \left( \dint_{{\varphi (\tau
,x}_{0}{)}}^{t}\partial _{x}G\left(s,\Phi _{{(\tau ,x}_{0}{)}}(s)\right)\,\,\mathrm{d}%
s\right) } \\ 
& =u_{0}(1)\exp \left( -\dint_{t_{0}}^{\tau }-\partial
_{x}G\left(s,\Phi _{{(\tau ,x}_{0}{)}}(s)\right)\,\mbox{d}s\right) \\ 
& =\underset{(t,x)\in \Omega _{2}}{\underset{(t,x)\rightarrow (\tau ,x_{0})}{%
\lim }}u(t,x).%
\end{array}%
\end{equation*}%
Now, let $(\tau ,x_{0})$ be a point on the curve $\Phi _{(t_{0},b)}$%
. As $u_{0}(b)=0,$ $\psi $\ is continuous and $\psi (\tau ,x_{0})=b,$ we get
\begin{equation*}
\begin{array}{ll}
\underset{(t,x)\in \Omega _{2}}{\underset{(t,x)\rightarrow (\tau ,x_{0})}{%
\lim }}u(t,x) & \quad =u_{0}\bigl(\psi (\tau ,x_{0})\bigr)\exp \left( -
\dint_{t_{0}}^{\tau }\partial _{x}G\left(s,\Phi _{(t,x)}(s)\right)\,\mbox{d}s\right) \\ 
& =u_{0}(b)\exp \left( -\dint_{t_{0}}^{\tau }\partial
_{x}G\left(s,\Phi _{(t,x)}(s)\right)\,\mbox{d}s\right) \\ 
& \quad =\underset{(t,x)\in \Omega _{3}}{\underset{(t,x)\rightarrow (\tau
,x_{0})}{\lim }}u(t,x) \\ 
& \quad =0.%
\end{array}%
\end{equation*}
\item We prove here that $\partial _{t}u \in\mathcal{C}^0_{MG}([t_{0},T],L^{1}(]1,b[))$. 
\begin{enumerate}[1.] 
\item Let $(t,x)\in \Omega _{1}$.
\begin{enumerate}[a.] 
\item  If $(t,x)$ belongs to one of the open connected components delimited by $\cup _{i=1,...,n}(\Sigma _{i}\cup \widetilde{\Sigma}_{i})$, $%
\widetilde{\Sigma}_{i}=\Phi _{(t_{i},1)},$ $\Sigma _{i}=\{t_{i}\}\times\lbrack 1,b\rbrack$, then $\varphi $ and $\Phi (.,.,s)$ are of class $C^{1}$ at $(t,x).$ The characteristics starting from the points $(t_{i},1)$ and those of foot $(\varphi (t,x),1)$ do not intersect. But, may be, there exist points of discontinuity of $G$ between $\varphi (t,x)$ and $t$ (see Figure. \ref{figure2}). So the functions $s\rightarrow \Phi _{(t,x)}(s)$ and $%
s\rightarrow \partial _{x}G(s,\Phi _{(t,x)}(s))$ are of class $C^{1}$ piecewise on the segment $[\varphi (t,x),t]$. By deriving the Lebesgue integral, $u$ is differentiable in time at the point $(t,x)$ and 
\begin{eqnarray}\label{Dtu1}
\partial _{t}u(t,x) &=&\exp\left(-\int_{\varphi (t,x)}^{t}\partial _{x}G\left(s,\Phi_{(t,x)}(s)\right)\ud s\right)\times\Bigg[R^{\prime }\bigl(\varphi (t,x)\bigr)\partial _{t}\varphi (t,x)\notag\\&&\qquad+R\bigl(\varphi (t,x)\bigr)\bigl(\partial _{t}\varphi(t,x)\partial _{x}G\bigl(\varphi (t,x),1\bigr)-\partial _{x}G(t,x)\bigr) \notag\\&&\qquad\quad-R\bigl(\varphi (t,x)\bigr)\int_{\varphi (t,x)}^{t}\partial_{x}^{2}G\left(s,\Phi_{(t,x)} (s)\right)\partial _{t}\Phi_{(t,x)} (s)\ud s\Bigg].
\end{eqnarray}%
\item  If $t\notin \{t_{i},i=1,...,n\}$ and $\varphi (t,x)=t_{j}$
with $t_{j}<t$, $\varphi $ is not derivable in time at the point $(t,x).$ There is a finite number, at most equal to $n,$
of such points. The function $\partial _{t}u(t,.)$ is thus defined almost
everywhere on $[1,\Phi _{(t_{0},1)}(t)]$ (see Figure. \ref{figure3}).
\item  If $t=t_{i}$ for one $i\in \{1,2,...n\},$ $\varphi $ is not derivable in time over the entire curve $\{(x,t_{i}),x\in \lbrack 1,\Phi_{(t_{0},1)}(t)]$ and admits right and left derivatives in time.

Finally, the function $\partial _{t}u$ is defined for any $t\in \lbrack t_{0},T]\backslash \{t_{i},$ $i=1,...,n\}$ and almost everywhere in $x$ in $\Omega _{1}.$
\end{enumerate}
\item  Let $(t,x)\in \Omega _{2}$.
\begin{enumerate}[a.] 
\item If $t\notin \{t_{i},i=1,...,n\}$, then $\psi $ and $\Phi (.,.,s)$ are of class $C^{1}$ at the point $(t,x).$ The function $s\rightarrow \Phi_{(t,x)}(s)$ is of class $C^{1}$ piecewise on $[t_{0},T]$. By derivation, $u$ is differentiable in time at the point $(t,x)$ and 
\begin{eqnarray}\label{Dtu2}
\partial _{t}u(t,x)&=&\exp \left(-\dint_{t_{0}}^{t}\partial _{x}G\left(s,\Phi_{(t,x)} (s)\right)\ud s\right)\times\notag\\
&& \qquad\qquad\Bigg[ \partial _{t}\psi (t,x)u_{0}^{\prime }\bigl(\psi (t,x)\bigr)- u_{0}\bigl(\psi(t,x)\bigr)\partial _{x}G(t,x)\notag\\
&&\qquad\qquad\quad-u_{0}(\psi(t,x)\dint_{t_{0}}^{t}\partial _{x}^{2}G\left(z,\Phi_{(t,x)}(s)\right)\partial _{t}\Phi _{(t,x)}(s)\,\mathrm{d}s\Bigg].
\end{eqnarray}%
\item  If $t\in \{t_{i},$ $i=1,...,n\},$\thinspace\ $u$ admits a right
and a left derivative in time on $(t,x)$ just like $\psi $ and $\Phi $.
\end{enumerate}
\end{enumerate}
Hence, for all $t\notin \{t_{i},$ $i=1,...,n\}$, the
mapping: $x\rightarrow \partial _{t}u(t,x)$ is defined almost everywhere on $[1,b]$ and verifies 
\begin{eqnarray*}
\dint_{1}^{b}\left\vert \partial _{t}u(t,x)\right\vert \ud x &\leq &M_{0}\Bigg[\sup_{\Omega _{2}}\big(|S(t,x)|+|\partial _{t}\psi (t,x)|\big)\Vert u_{0}\Vert _{W^{1,1}}\\
&&C_{1}\sup_{\Omega _{1}}\big(|F(t,x)|+|\partial _{t}\varphi (t,x)|\big)
\bigg(\Vert \beta \Vert _{\infty }\Vert w\Vert _{E_{MG}}+\Vert f\Vert _{-\mathcal{C}_{MG}^{1}}\bigg)\Bigg],
\end{eqnarray*}%
where 
\begin{eqnarray}\label{F}
F(t,x) &=&-\partial _{x}G(t,x)+\partial _{t}\varphi (t,x)\partial
_{x}G\bigl(\varphi (t,x),1\bigr)\notag\\
&&\qquad\qquad\qquad\qquad\qquad-\int_{\varphi (t,x)}^{t}\partial _{x}^{2}G\left(s,\Phi_{(t,x)}
(s)\right)\partial _{t}\Phi_{(t,x)} (s)\ud s,  \notag
\end{eqnarray}%
and 
\begin{eqnarray}\label{S}
S(t,x) &=&\partial _{x}G(t,x)+\int_{t_{0}}^{t}\partial _{x}^{2}G\left(z,\Phi
_{(t,x)}(s)\right)\partial _{t}\Phi _{(t,x)}(s)\,\,\mathrm{d}s.
\end{eqnarray}%
Thus $\partial _{t}u(t,.)\in L^{1}(]1,b[).$ 

In addition, the mapping: $t\rightarrow\partial_{t}u(t,.)$ is continuous almost everywhere from $[t_{0},T]$ in $%
L^{1}(]1,b[).$ Indeed, for any $t\in\lbrack t_{0},T]\backslash\cup_{i=1}^{n}\{t_{i}\}$ and any sequence $(\hat{t}_{n})_{n\in \mathbb{N}%
}\subset \lbrack t_{0},T]\backslash \cup _{i=1}^{n}\{t_{i}\}$ converging to $%
t$, we have
\begin{eqnarray*}
\left\Vert \partial _{t}u(\hat{t}_{n},.)-\partial _{t}u(t,\,.)\right\Vert
_{L^{1}} &\leq&M_{0}C_{1}\Bigg[\sup_{\Omega _{1}}|F(t,x)|\bigg(\Vert\beta \Vert _{\infty}\left\Vert w(\hat{t}_{n},.)-w(t,.)\right\Vert_{L^{1}}\\
&&+\left\vert f(\hat{t}_{n})-f(t)\right\vert \bigg)+\sup_{\Omega _{1}}|\partial _{t}\varphi (t,x)|\bigg(\left\vert f^{\prime }(\hat{t}_{n})-f^{\prime}(t)\right\vert\\
&&\qquad\qquad\qquad+\Vert \beta \Vert_{\infty }\left\Vert \partial _{t}w(\hat{t}_{n},.)-\partial_{t}w(t,.)\right\Vert _{L^{1}}\bigg)\Bigg].
\end{eqnarray*}
We conclude by using the regularity of $f$ and $w$.
\end{enumerate}
\end{proof}
\begin{lemma}
Under the hypothesis of Theorem \ref{th_existence_unicite}, the operator $\mathcal{T}$ admits a unique fixed point.
\end{lemma}
\begin{proof}
\begin{enumerate}[Step 1.]
\item Let $w_{1},\,w_{2}$ be given in $E(u_{0})$. We have 
\begin{equation*}
\mathcal{T}(w_{1})-\mathcal{T}(w_{2})=0\text{ on }\Omega _{2}\cup \Omega _{3},
\end{equation*}%
and then for all $(t,x)\in \mathcal{Q}$ we have
\begin{eqnarray*}
|\left( \mathcal{T}(w_{1})-\mathcal{T}(w_{2})\right) (t,x)|&\leq&\dfrac{M_{0}C_{1}}{b-1}\Vert \beta \Vert _{\infty } \Vert w_{1}\bigl(\varphi(t,x),.\bigr)-w_{2}\bigl(\varphi (t,x),.\bigr)\Vert _{L^{1}}\mathds{1}_{\Omega_{1}}(t,x),
\end{eqnarray*}%
where $M_{0}$ and $C_{1}$ are given by (\ref{M0}). It follows that 
\begin{eqnarray*}
\Vert \left( \mathcal{T}(w_{1})-\mathcal{T}(w_{2})\right) (t,.)\Vert _{L^{1}}
& \leq& \dfrac{M_{0}C_{1}}{b-1}\Vert \beta \Vert _{\infty }\displaystyle{\int_{1}^{\Phi_{(t_{0},1)}(t)}\Vert \left( w_{1}-w_{2}\right) \bigl(\varphi (t,x),.\bigr)\Vert_{L^{1}}\ud x} \\ 
&\leq &M_{0}C_{1}\Vert \beta \Vert _{\infty }\Vert w_{1}-w_{2}\Vert _{E(u_{0})}.%
\end{eqnarray*}%
On another hand, a technical calculation gives for almost all $(t,x)$
\begin{eqnarray*}
\partial _{t}\left( \mathcal{T}(w_{1})-\mathcal{T}(w_{2})\right) (t,x) &=&\mathds{1}_{\Omega_{1}}(t,x)\times\exp \left( -\displaystyle{\int_{\varphi (t,x)}^{t}\partial _{x}G\left(s,\Phi_{(t,x)}(s)\right)\ud s}\right) \times\\
&&\quad\bigg[\dfrac{\partial _{t}\varphi (t,x)C_1}{b-1}\displaystyle{\int_{1}^{b}\beta (y)\partial_{t}(w_{1}-w_{2})\bigl(\varphi (t,x),y\bigr)\ud y} \\ 
&& \qquad+\dfrac{F(t,x)C_1}{b-1}\displaystyle{\int_{1}^{b}\beta(y)(w_{1}-w_{2})\bigl(\varphi (t,x),y\bigr)\,\mbox{d}y}\bigg]
 , 
\end{eqnarray*}%
where $F$ is defined by (\ref{F}). It follows that 
\begin{eqnarray*}
\Vert \partial _{t}\left( \mathcal{T}(w_{1})-\mathcal{T}(w_{2})\right)
(t,.)\Vert _{L^{1}}&\leq& \dfrac{M_{0}C_{1}\left( C_{2}+C_{3}\right)}{b-1}\Vert \beta \Vert _{\infty }\Phi _{(t_{0},1)}(t)\Vert w_{1}-w_{2}\Vert _{E(u_{0})},
\end{eqnarray*}%
with $C_{2}=\sup_{Q}\{1,\,|\partial _{t}\varphi (t,x)|\}$ and $%
C_{3}=\sup_{Q}|F(t,x)|$.\newline
For $t^{\prime }$ such that $\Phi _{(t_{0},1)}(t^{\prime })=\dfrac{b-1}{%
2M_{0}C_{1}\Vert \beta \Vert _{\infty }\left( 1+C_{2}+C_{3}\right) }+1,$ we get
\begin{eqnarray*}
\Vert \mathcal{T}(w_{1})-\mathcal{T}(w_{2})\Vert _{E(u_{0})}&\leq&
C\Vert \beta \Vert _{\infty } \left(\Phi _{(t_{0},1)}(t^{\prime })-1\right) \Vert w_{1}-w_{2}\Vert _{E(u_{0})},
\end{eqnarray*}%
Where $$C=\dfrac{M_{0}C_{1}\left( 1+C_{2}+C_{3}\right)}{b-1}.$$
Let $\tilde{t}=\Phi _{(t_{0},1)}^{-1}\left( \dfrac{b-1}{2M_{0}C_{1}\Vert \beta
\Vert _{\infty }\left( 1+C_{2}+C_{3}\right) }+1\right) $ and $t^{\ast }=\min(\tilde{t},T)$. The operator $%
\mathcal{T}$ is a contracting operator on 
\begin{equation*}
E_{t^{\ast }}(u_{0})=\big\lbrace w\in \mathcal{C}_{MG}^{1}([t_{0},t^{\star
}],\,L^{1}([1,b]))\ ; w(t_{0},\,.)=u_{0}\big\rbrace.
\end{equation*}
So it admits a unique fixed point.
\item The contraction constant, in the step 1, do not depend on $u_{0}$, we
can therefore repeat the process of the first step with an initial condition 
$u_{0}^{\ast }=u(t^{\ast },\,.)$ and the compatibility condition at $t^{\ast
}$. We obtain a new point $t^{2\ast }$ such that the operator $\mathcal{T}$
is contracting on 
\begin{equation*}
E_{t^{2\ast }}(u_{0}^{\ast })=\big\lbrace w\in \mathcal{C}%
_{MG}^{1}([t^{\star },t^{2\star }],\,L^{1}([1,b]))\;w(t^{\star
},\,.)=u(t^{\star },\,.)\big\rbrace,
\end{equation*}%
so it admits a unique fixed point on $[t^{\ast },\,t^{2\ast }]$. We iterate the process until we cover $[t_{0},T]$.
\end{enumerate}
\end{proof}
\begin{lemma}
Under the hypothesis of Theorem \ref{th_existence_unicite}, the fixed point $u$ of the two previous lemma is solution of the problem (\ref{Eq1})-(\ref{Eq6}).
\end{lemma}
\begin{proof}
By the construction (\ref{Cons}) and by taking into account the respective regularities, we have\\ 
$\begin{array}{llll}
\partial_x(Gu)(t,x)&= & 0,\quad  \text{in }\Omega_3,\\
\partial_x(Gu)(t,x)&=&\exp\left(-\dint_{\phi(t,x)}^t\partial_x G\left(s,\Phi_{(t,x)}(s)\right)\ud s\right)\times \\
&&\Biggl[R\bigl(\phi(t,x)\bigr)G(t,x)\partial_x\phi(t,x)\partial_xG\bigl(\phi(t,x),1\bigr)\\
&&-R\bigl(\phi(t,x)\bigr)G(t,x)\dint_{\phi(t,x)}^t\partial_x^2G\left(s,\Phi_{(t,x)}(s)\right)\partial_x\Phi_{(t,x)}(s)\ud s\\
&&+G(t,x)\partial_x \phi(t,x)R^{\prime}\bigl(\phi(t,x)\bigr)-R\bigl(\phi(t,x)\bigr)\partial_xG(t,x)\Biggr],\quad\text{in } \Omega_1,\\
\partial_x(Gu)(t,x)&=&\exp\left(-\dint_{t_0}^t \partial_xG\left(s,\Phi_{(t,x)}(s)\right) \ud s\right)\times\\
&&\Biggl[G(t,x)\partial_x\psi(t,x)u'_0\bigl(\psi(t,x)\bigr)-u_0\bigl(\psi(t,x)\bigr)\partial_x G(t,x)\\
&&- u_0\bigl(\psi(t,x)\bigr)G(t,x) \dint_{t_0}^t\partial_x^2 G\left(s,\Phi_{(t,x)}(s)\right)\partial_t\Phi_{(t,x)}(s) \ud s\Biggr],\quad\text{in } \Omega_2.
\end{array}$\\
Using the derivative in time (\ref{Dtu1})-(\ref{Dtu2}), we get 
\begin{equation}
\partial _{t}u+\partial _{x}(Gu)=\left\{ 
\begin{array}{ll}
R^{\prime }\bigl(\varphi (t,x)\bigr)\exp \left( -\dint_{\varphi (t,x)}^{t}\partial _{x}G\left(s,\Phi _{(t,x)}(s)\right)\ud s\right)\times &\\
\qquad\qquad\qquad\qquad\bigg[\partial_{t}\varphi (t,x)+G(t,x)\partial _{x}\varphi (t,x)\bigg], & (t,x)\in \Omega _{1}\\ 
u_{0}\bigl(\psi (t,x)\bigr)\exp \left( -\dint_{t_{0}}^{t}\partial_{x}G\left(s,\Phi _{(t,x)}(s)\right)\ud s\right)\times &\\
\qquad\qquad\qquad\qquad\bigg[\partial _{t}\psi(t,x)+G(t,x)\partial _{x}\psi (t,x)\bigg], & (t,x)\in \Omega _{2} \\ 
0. & (t,x)\in \Omega _{3}%
\end{array}%
\right.  \label{EEE}
\end{equation}%
Moreover, since the functions $\varphi $ and $\psi $ are constant along the characteristics curves, we have
\begin{equation}
\frac{d}{ds}\varphi \bigl(s,\Phi (s;t,x)\bigr)=0\text{ \ \ \ \ and \ \ }\frac{d}{ds}\psi \bigl(s,\Phi (s;t,x)\bigr)=0\label{EE1}.
\end{equation}%
Then, we get
\begin{equation}
\partial _{s}\varphi \bigl(s,\Phi (s;t,x)\bigr)+G\bigl(t,\Phi (s;t,x)\bigr)\partial _{x}\varphi\bigl(s,\Phi (s;t,x)\bigr) =0  \label{EE2} \\
\end{equation}
and
\begin{equation*}
\partial _{s}\psi \bigl(s,\Phi (s;t,x)\bigr)+G\bigl(t,\Phi (s;t,x)\bigr)\partial _{x}\psi\bigl(s,\Phi (s;t,x)\bigr) =0. 
\end{equation*}
Using (\ref{EEE}), (\ref{EE1}) and (\ref{EE2}), we deduce that
\begin{equation}
\partial _{t}u(t,x)+\partial _{x}(G(t,x)u(t,x))=0.
\end{equation}%
Finally, since $\partial _{t}u\in L^{1}(]1,b[)$, we get $\partial
_{x}(Gu)=-\partial _{t}u\in \mathcal{C}([t_{0},t_{1}],L^{1}(]1,b[))$ and $u$
is a solution of our problem. 
\end{proof}
\begin{lemma}
Under hypothesis of Theorem \ref{th_existence_unicite}, the problem (\ref{Eq1})--(\ref{Eq6}) has a unique solution.
\end{lemma}
\begin{proof}
Let $u$ be a solution of (\ref{Eq1})--(\ref{Eq6}). We decompose $u$ in the form
\begin{equation}
u=u_{1}+u_{2},  \label{decomp}
\end{equation}%
where $u_{1}$ and $u_{2}$ are solutions of (\ref{Eq1})--(\ref{Eq3}) with respectively $%
f=f_{1}=0$, $u_{1}(0,.)=u_{0}$ and $f_{2}=f$, $u_{2}(0,.)=0$. 
For $n\in \mathbb{N}^{\ast }$, we introduce the function $\Gamma _{n}$ defined on $%
\mathbb{R}$ by 
\begin{equation*}
\Gamma _{n}(x)=\left\{ 
\begin{array}{lc}
\,\,\,\,1 & \text{if}\,\,x\geq \frac{1}{n}, \\ 
-1, & \text{if}\,\,x\leq -\frac{1}{n} \\ 
\,\,\,\,x & \text{if}\,\,-\frac{1}{n}\leq x\leq \frac{1}{n}.%
\end{array}%
\right.
\end{equation*}%
Multiplying $\partial _{t}u_{i}+\partial _{x}\left[ Gu_{i}\right] $, for $%
i=1,2$, by $\Gamma _{n}(u_{i}(t,x))$ and integrating on $[1,b]$, we get
\begin{equation*}
\dint_{1}^{b}\partial _{t}u_{i}(t,x)\Gamma
_{n}\bigl(u_{i}(t,x)\bigr)+\partial _{x}\left[ G(t,x)u_{i}(t,x)\right] \Gamma
_{n}\bigl(u_{i}(t,x)\bigr)\ud x=0.
\end{equation*}%
Integrating by parts, making $n$ tend towards infinity and applying the
dominated convergence theorem, we get for $i=1,2$ 
\begin{equation*}
\begin{array}{ll}
\partial _{t}||u_{i}(t,.)||_{L^{1}(1,b)} & \leq \bigl|\displaystyle{%
\dint_{1}^{b}\beta (x)u_{i}(t,x)\,\mbox{d}x}\bigr|+\bigl|f_{i}(t)\bigr| \\ 
& \leq ||\beta ||_{\infty }||u_{i}(t,.)||_{L^{1}(1,b)}+\bigl|f_{i}(t)\bigr|%
\end{array}.%
\end{equation*}%
Integrating with respect to $t$, for all $i=1,2$, we obtain
\begin{equation*}
\Vert u_{i}(t,.)\Vert _{L^{1}(1,b)}\leq \Vert \beta \Vert _{\infty }%
\displaystyle{\int_{t_{0}}^{t}\Vert u_{i}(s,.)\Vert _{L^{1}(1,b)}\,\mbox{d}%
s+\int_{t_{0}}^{t}\bigl|f_{i}(s)\bigr|\,\mbox{d}s}+\Vert u_{i}(0,x)\Vert
_{L^{1}(1,b)}.
\end{equation*}%
Then by applying the Gronwall's Lemma, for $i=1,2$, we get
\begin{equation}
\Vert u_{i}(t,.)\Vert _{L^{1}(1,b)}\leq \mbox{e}^{\Vert \beta \Vert _{\infty
}(T-t_{0})}\bigl(||f_{i}||_{L^{1}(t_{0},T)}+||u_{i}(0,x)||_{L^{1}(1,b)}\bigr).
\label{estimation1}
\end{equation}
We denote by $\mathcal{T}_{u_{0},f}$ the operator constructed in Lemma \ref{const}. Hence for $u_{0}\in W^{1,1}(]1,b[)$ verifying the compatibility condition $(Gu_{0})(t_{0},1)=\int_{1}^{b}\beta (x)u_{0}(x)dx,$ the operator $%
\mathcal{T}_{1}=\mathcal{T}_{u_{0},0}:u_{0}\mapsto u_{1}$ satisfies 
\begin{equation}
||\mathcal{T}_{u_{0},0}u_{0}||_{\mathcal{C}([t_{0},T],L^{1}(1,b))}\leq %
\mbox{e}^{||\beta ||_{\infty }(T-t_{0})}||u_{0}||_{L^{1}(1,b)}.  \label{est1}
\end{equation}%
Now, for $f\in \mathcal{C}_{MG}^{1}([t_{0},T])$ such that $f(t_{0})=0,$ the
function still denoted by $f$ defined by $f(t,x)=f(t)\ $, belongs to $E(u_{0}=0)%
\mathcal{\ }$and the operator $\mathcal{T}_{2}$ defined by $\mathcal{T}_{2}=%
\mathcal{T}_{0,f}:f\mapsto u_{2}$ verifies
\begin{equation}
||\mathcal{T}_{2}f||_{\mathcal{C}([t_{0},T],L^{1}(1,b))}\leq \mbox{e}%
^{||\beta ||_{\infty }(T-t_{0})}||f||_{L^{1}(t_{0},T)}.  \label{est2}
\end{equation}%
Finally, if there exist two solutions $v_{1}$ and $v_{2}$ of (\ref{Eq1})--(\ref{Eq3}), then $%
u=v_{1}-v_{2}$ is a solution with initial data $u_{0}$ and second member $f$%
\ null. We decompose $u$ as in (\ref{decomp}) and we deduce from (\ref{est1}, %
\ref{est2}) that $u_{1}=u_{2}=0$ and the uniqueness follows.
\end{proof}
\section{Numerical analysis}
In this section, we will present some numerical tests that would show the interest of drugs on tumor growth. With $G$ given by (\ref{G}), it is not possible to compute an exact solution for the system (\ref{Eq1}-\ref{Eq6}) except in very special cases. So, it is necessary to implement an approximation strategy. We use the method of characteristics to solve the transport equation, a 4th order Runge Kutta scheme for the ODE and a trapezoidal quadrature formula to approximate the integrals. To follow the evolution of the disease, we introduce as in \cite{Barbolosi}, the metastatic index defined by 
\begin{equation}
MI_{b_{\min }}(t)=\int_{b_{\min }}^{b}u(t,x)\ud x  \label{MI}.
\end{equation}%
It evaluates the number of metastases at different times and for different minimal sizes $b_{\min }$. If $b_{\min }=1,$ all the metastases present are measured. For $b_{\min }\geq 10^{8}$, only detectable tumors by medical imaging are considered \cite{Barbolosi}. To measure the influence of drugs, we will present the evolution of this index without treatment, with chemotherapy, and with mixed treatments.

\subsection{Method of characteristics}
\subsubsection{Time discretization}
We consider a time interval $[0,T]$ and an interval $[t_{0},\overline{T}]$ strictly included in $[0,T]$, where $t_{0}$ and $\overline{T}$ are respectively the start and the end of treatments. We notice that we are in the presence of large variations in time scales: the radiation time is of the order of minutes, the chemo time is counted in hours, the treatment lasts several months and the disease is monitored for a few years. We then opt for a time discretization with variable steps. Let $(t_{n})_{0\leq n\leq N}$ be a nonuniform subdivision of $[0,T]$ with steps $k_{n}=t_{n+1}-t_{n}$. Care will be taken to consider the points of discontinuity of $G$ as nodes of the mesh.  

\subsubsection{Space discretization}
The spatial discretization is based on the resolution of the differential equation
\begin{equation}
\left\{ 
\begin{array}{lcl}
\frac{d}{dt}\Phi(t)&=&G(t,\Phi(t)), \\ 
\Phi (s)&=&x.%
\end{array}%
\right.  \label{CA1}
\end{equation}%
Along this curve, the solution of the EDP is given by%
\begin{equation}
u(t,\Phi _{(s,x)}(t))=u(s,x)\exp \Bigg(-\displaystyle{\int_{s}^{t}\partial
_{x}G(\tau ,\Phi _{(s,x)}(\tau ))\,\mbox{d}\tau }\Bigg).  \label{INT1}
\end{equation}
To simplify the notations, we introduce the function $y_{p}(t)$ defined by%
\begin{equation*}
y_{p}(t)=\left\{ 
\begin{array}{ll}
x_{p}(t)=b^{1-{\mbox{e}}^{-at}} & \mbox{if}\,\,t\in \lbrack 0,\,t_{0}], \\ 
\Phi _{(t_{0},1)}(t) & \text{if not}.
\end{array}%
\right.
\end{equation*}%
We define a spacial mesh $(x_{i}^{n})_{i\geq 1}$ at time $t_{n}$ by 
\begin{eqnarray}
x_{i}^{n} &=&\Phi _{(t_{n-1},x_{i-1}^{n-1})}(t_{n}),\text{ \ }\forall i\geq
2,\text{ \ }n\in \lbrack \lbrack 2,N+1]],  \label{un} \\
x_{1}^{n} &=&1,\text{ \ }\forall \text{ }n,  \label{deux} \\
x_{i}^{0} &=&y_{p}(t_{i}),\text{ \ }\forall i.  \label{quatre}
\end{eqnarray}%
Equality (\ref{deux}) means that $(t,1)$ is the entry point for all the characteristic curves in
\begin{equation*}
\Omega _{1}=\big\lbrace(t,x)\in \mathcal{Q}\,\,\text{\textit{/ }}\,x<\Phi
_{(t_{0},1)}(t)\big\rbrace.
\end{equation*}\\ 
We denote the step of the space discretization by 
\begin{equation*}
h_{i}^{n}=x_{i}^{n}-x_{i-1}^{n}.
\end{equation*}%
Using a trapezoidal quadrature formula in the expression (\ref{INT1}), we get %
\begin{equation}
\left\{ 
\begin{array}{ll}
u_{1}^{1}=\dfrac{\beta (y_{p}(0))}{G(1,1)}, &  \\ 
u_{i}^{n}=0, & i\geq n,\\ 
u_{i}^{n}=u_{i-1}^{n-1} {\mbox{e}}^{-%
\mathfrak{S}\left( t_{n-1},t_{n},x_{i-1}^{n-1},x_{i}^{n},\partial
_{x}G\right)}, & n=2,...,N,i=2,...,n, \\ 
u_{1}^{n}=\dfrac{1}{G(t_{n},1)}\Bigg(\beta (y_{p}(t_{n}))+\mathfrak{Q}\left(
(h_{i}^{n})_{i\geq 2},(u_{i}^{n})_{i\geq 2}\right) \Bigg), & n=2,...,N,%
\end{array}%
\right.  \label{DISCRET}
\end{equation}%
where 
\begin{eqnarray*}
\mathfrak{S}\left( t_{n-1},t_{n},x_{i-1}^{n-1},x_{i}^{n},\partial
_{x}G\right) &=&{\frac{1}{2}}k_{n}\bigg(\partial
_{x}G(t_{n},x_{i}^{n})+\partial _{x}G(t_{n-1},x_{i-1}^{n-1})\bigg), \\
&&\text{and} \\
\mathfrak{Q}\left( (h_{i}^{n})_{i\geq 2}, (u_{i}^{n})_{i\geq 2} \right) &=&h_{2}^{n}\beta
(x_{2}^{n})u_{2}^{n}+\sum_{i=3}^{n}\left( {\frac{1}{2}}h_{i}^{n}(\beta
(x_{i}^{n})u_{i}^{n}+\beta (x_{i-1}^{n})u_{i-1}^{n}\right),
\end{eqnarray*}%
are respectively the quadrature formula in $
\dint_{t_{n-1}}^{t_{n}}\partial _{x}G(z,\Phi _{(t_{n-1},x_{i-1}^{n-1})}(z))\,%
\ud z$ and in the nonlocal boundary condition.
The second equation of (\ref{DISCRET}) means that the solution vanishes on
\begin{equation*}
\Omega _{3}=\big\lbrace(t,x)\in \mathcal{Q}\,\,\text{\textit{/ }}\,x>\Phi
_{(t_{0},b)}(t)\big\rbrace.
\end{equation*}
\subsection{ Numerical results}
It is worthy to note that all our calculations are carried out with MATLAB software. The time step varies between 3 min and 1 hour. The spatial variable is a matrix of the order of $2\times 10^5$ with coefficients varying between $1$ and $10^{12}$, depending on whether one is in the regions with or without treatment. 
\subsubsection{Validation tests}
\begin{figure}[h!]
    \centering
    \includegraphics[width=0.5\linewidth,height=0.2\linewidth]{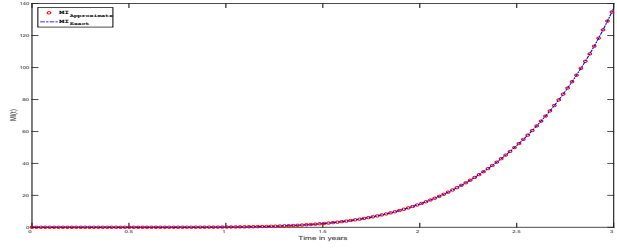}
    \caption{Exact and approximate MI without treatment: $a=0.00286$, $b=10^{12}$, $b_{min}=1$, $m=5.3\times 10^{-8}$, $\alpha=\frac{2}{3}$, $k_{n}=1$ hour for all $n$.}
  \label{figWT}
\end{figure}

\begin{figure}[h!] 
\centering
\begin{subfigure}[b]{.46\linewidth}
   \centering
    \includegraphics[width=0.85\linewidth,height=0.45\linewidth]{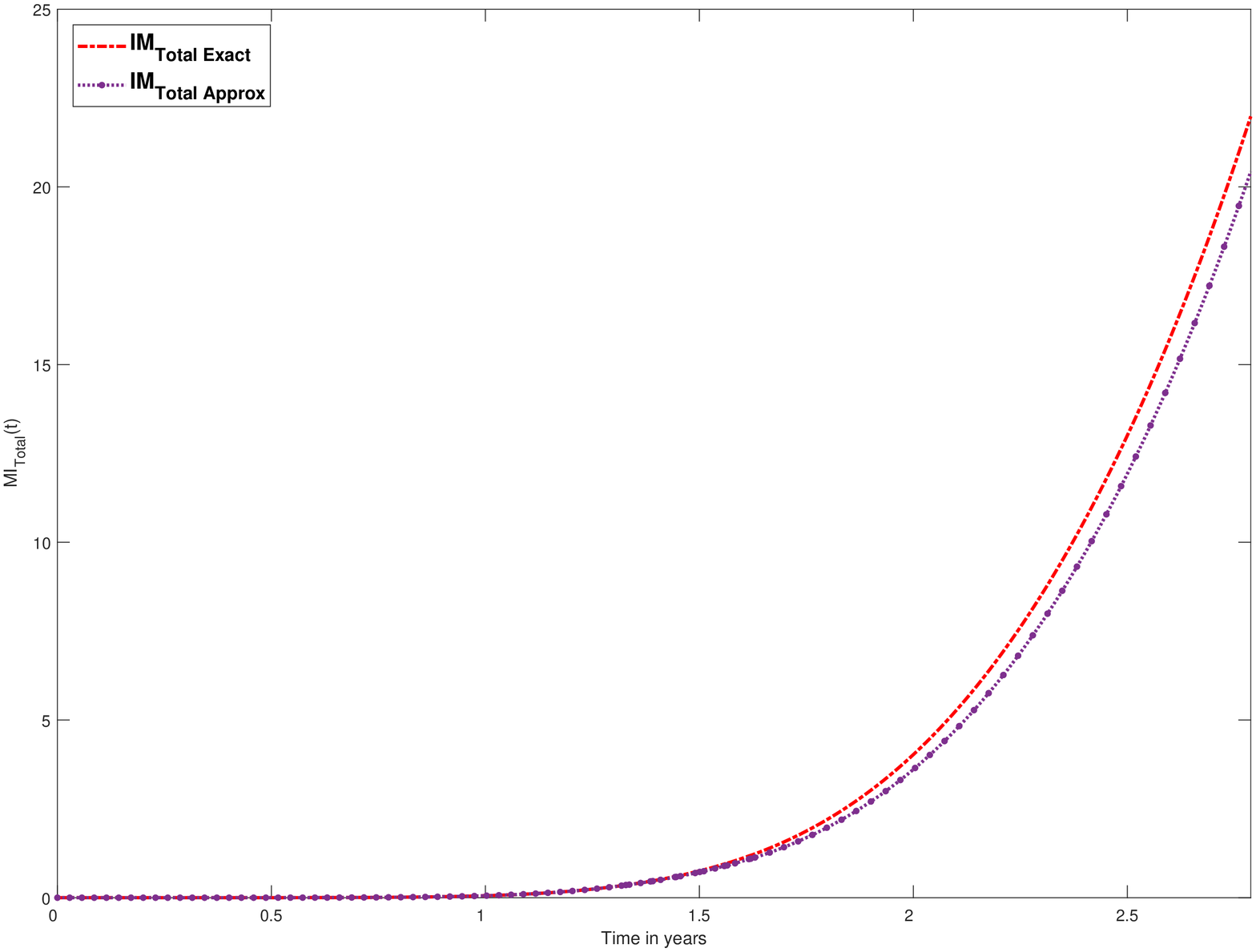}
    \caption{Total MI ($b_{min}=1$),}
  \label{fig5}
\end{subfigure}
\hfill
\begin{subfigure}[b]{.46\linewidth}
   \centering
    \includegraphics[width=0.85\linewidth,height=0.45\linewidth]{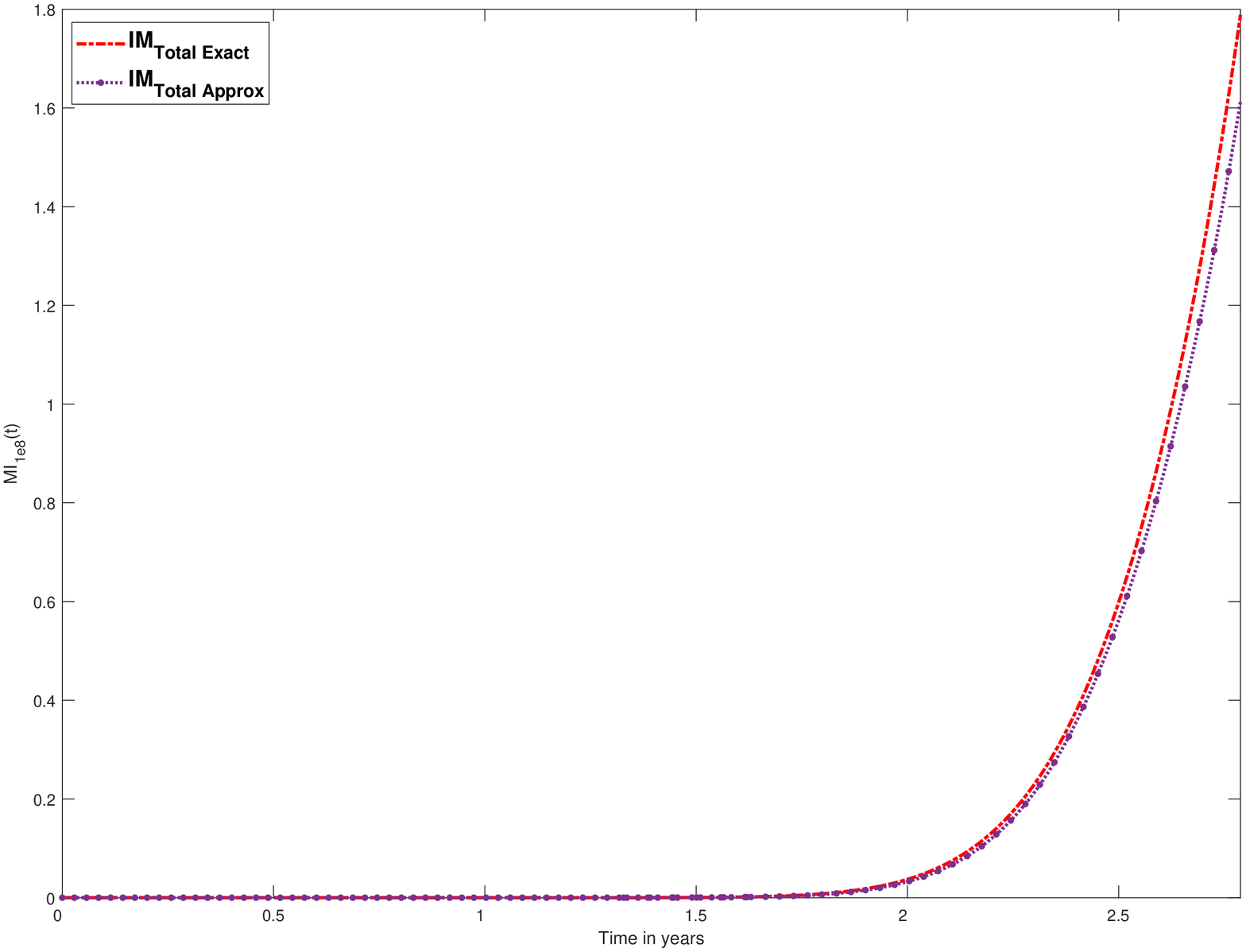}
    \caption{Visible MI ($b_{min}=10^8$),}
  \label{fig6}
\end{subfigure}
\caption{Comparison of exact and approximate metastatic index with Chemotherapy.}
\end{figure}

\begin{figure}[h!] 
    \centering
    \includegraphics[width=0.5\linewidth,height=0.2\linewidth]{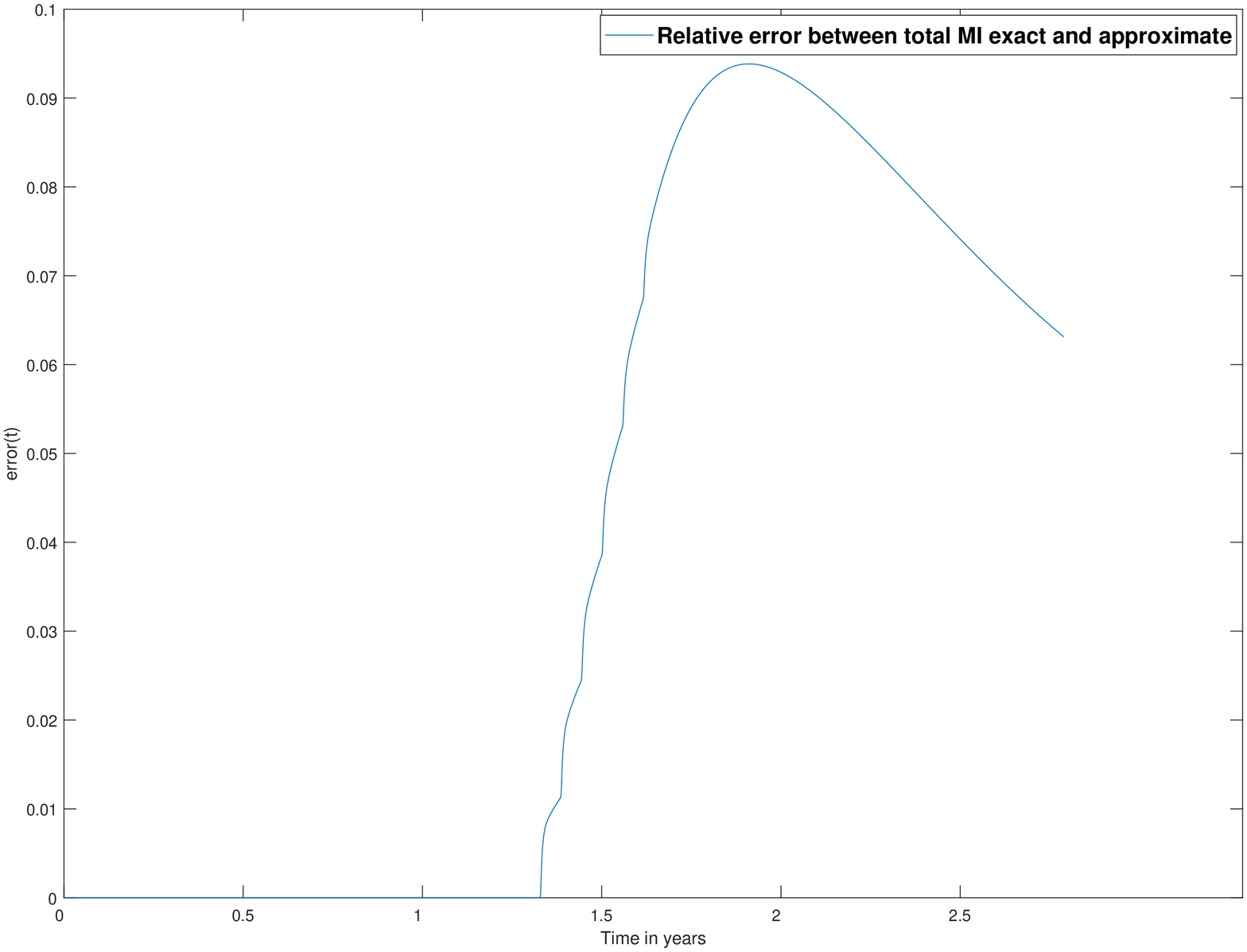}
    \caption{Relative error between exact and approximate total MI with chemotherapy}
  \label{err2}
 \end{figure}

\begin{figure}[h!]
\begin{subfigure}[b]{.49\linewidth}
    \centering
    \includegraphics[width=0.85\linewidth,height=0.45\linewidth]{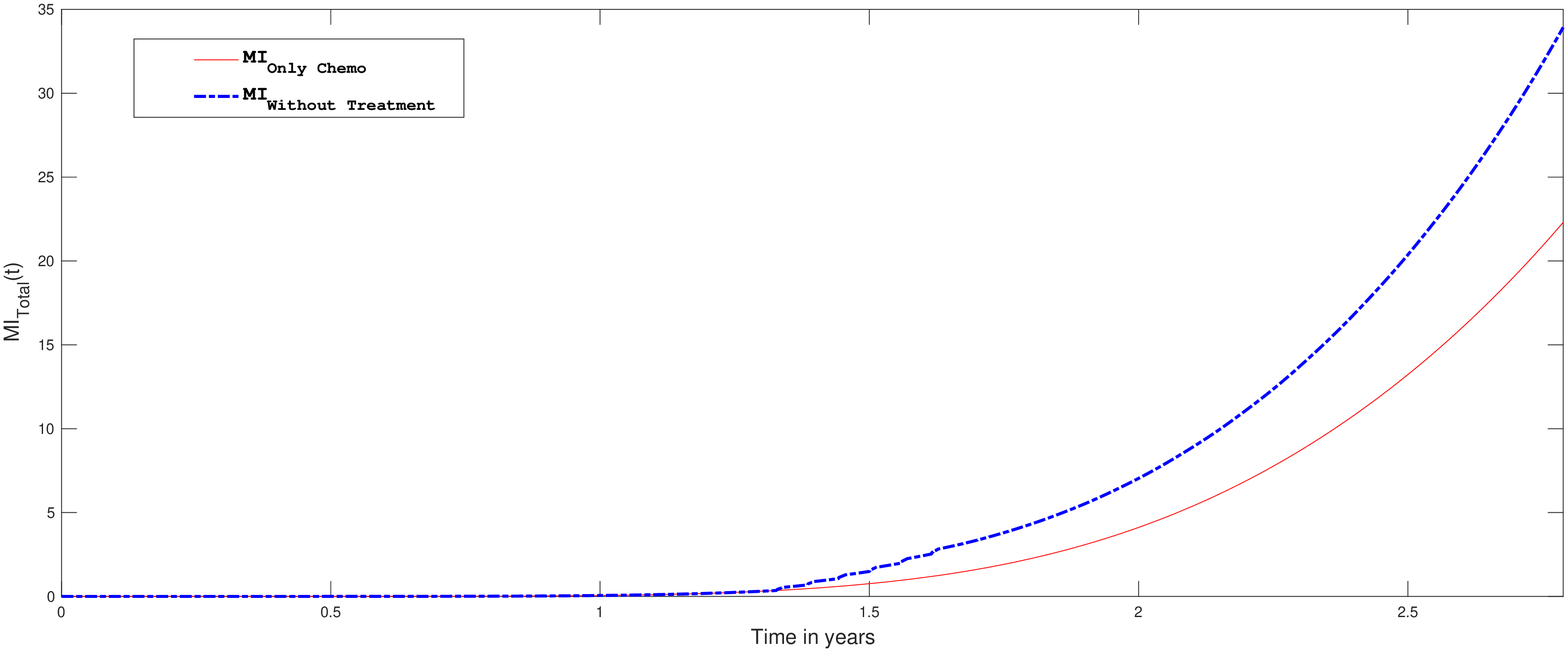}
    \caption{Total MI,}
  \label{fig7}
\end{subfigure}
\hfill
\begin{subfigure}[b]{.49\linewidth}
    \centering
    \includegraphics[width=0.85\linewidth,height=0.45\linewidth]{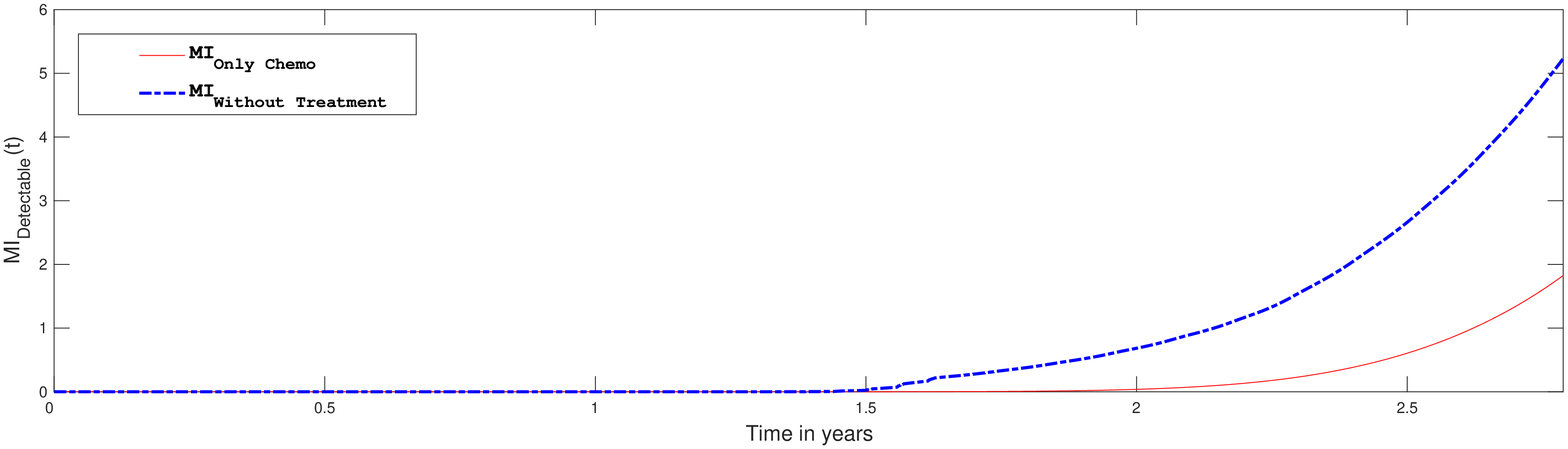}
    \caption{Detectable MI,}
  \label{fig8}
\end{subfigure}
\caption{Metastatic indexes with and without chemotherapy.}
\end{figure}

In order to measure the performance of our numerical strategy, we consider two tests with known analytical solutions.
\paragraph{First test : Without treatment}
The characteristics curves are solutions of (\ref{Eq1}) and are given by  
\begin{equation}
\Phi(t,s,x)=b\left(\dfrac{x}{b}\right)^{\e^{-a(t-s)}}.
\end{equation}
We present in Figure \ref{figWT} the exact and approximate metastatic index and observe that the numerical solution is in perfect agreement with the exact one.
\paragraph{Second test : Only chemotherapy}
We consider the case with with
\begin{equation*}
K(x)=x,\quad \Gamma(t)=\gamma C(t)\quad\text{and}\quad C(t)=c_{1}(t),
\end{equation*}%
where $c_{1}(t)$ is the drug's concentration defined by (\ref{C1}).
The exact solution for the ODE (\ref{Eq6}) is 
\begin{equation}
\Phi(t,s,x)=b\left(\dfrac{x}{b}\right)^{\e^{-a(t-s)}}\e^{-\gamma\e^{-at}\left(I(t)-I(s)\right)},
\end{equation}
with $I(t)=\int_0^t C(\tau)\e^{a\tau}\ud \tau$. Using the following parameters 
\begin{equation*}
    a=0.00286,\quad b=10^{12},\quad m=5.3\times 10^{-8},\quad \alpha=0.55,
\end{equation*}
we present in Figure \ref{fig5} and Figure \ref{fig6}, the approximation error on the metastatic index MI, respectively for detectable and total tumors. We can see the performance of our methodology despite the disparity in the parameters and the size of the matrix. We also see that MI is larger with ${b_{min} = 1}$, which highlights the existence of not detectable metastases by imaging. Figure \ref{err2} shows that the relative error is of the order of $9\%$ and decreases after the end of treatment i.e. when $G$ becomes regular. In Figure \ref{fig7} and Figure \ref{fig8}, we highlight the positive effect of the treatments on the metastatic indexes. The difference is of order of $27\%$ for all metastasis and $33\%$ for detectable metastasis. All the results of this test are consistent with those of \cite{FV} and confirm the validation of our calculation code.

\hfill

\subsubsection{Significant tests}
\begin{figure}[h!]
\begin{subfigure}[b]{.49\linewidth}
    \centering
    \includegraphics[width=0.85\linewidth,height=0.45\linewidth]{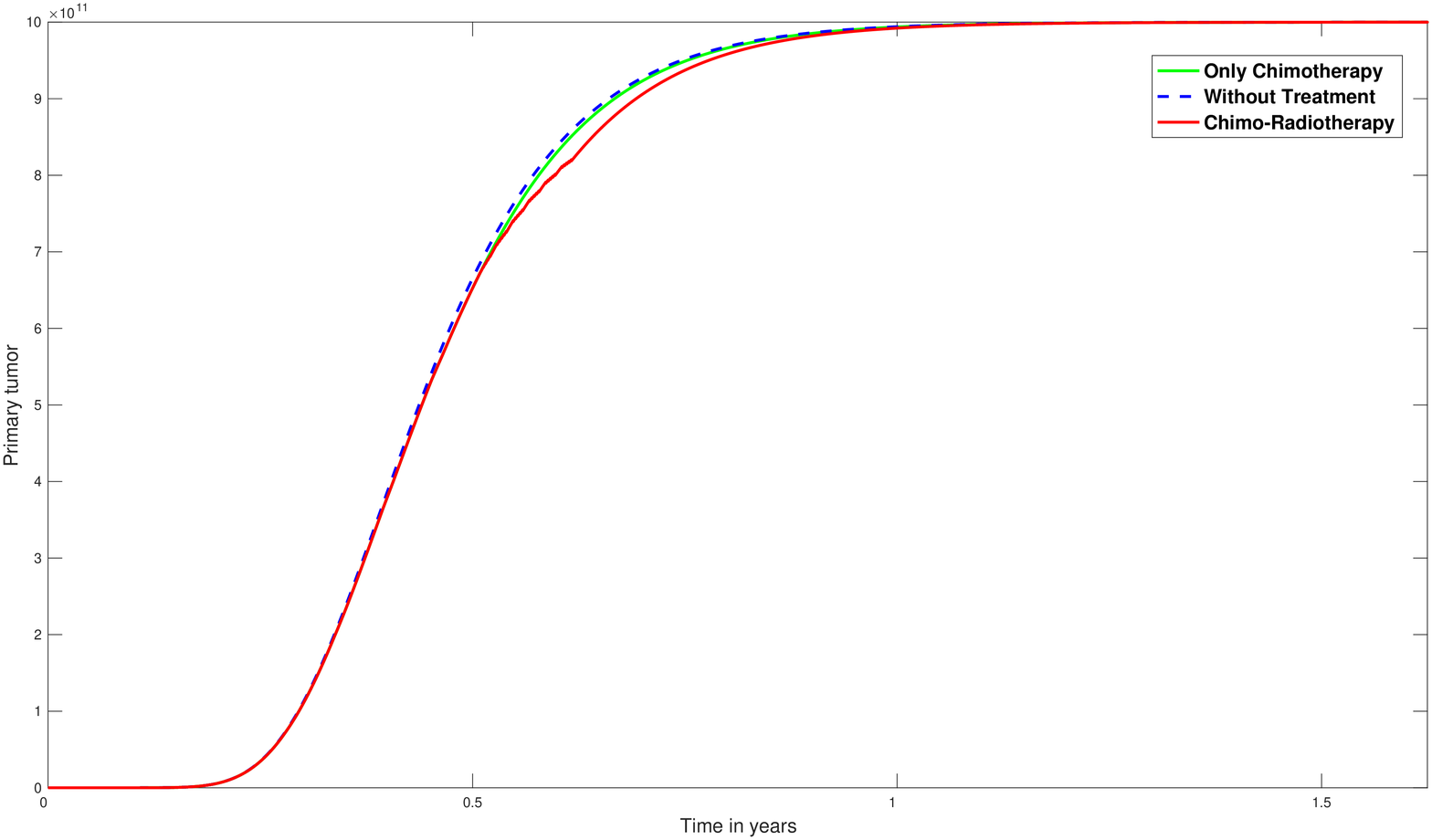}
    \caption{Primary tumor,}
  \label{fig9}
\end{subfigure}
\hfill
\begin{subfigure}[b]{.49\linewidth}
    \centering
    \includegraphics[width=0.85\linewidth,height=0.45\linewidth]{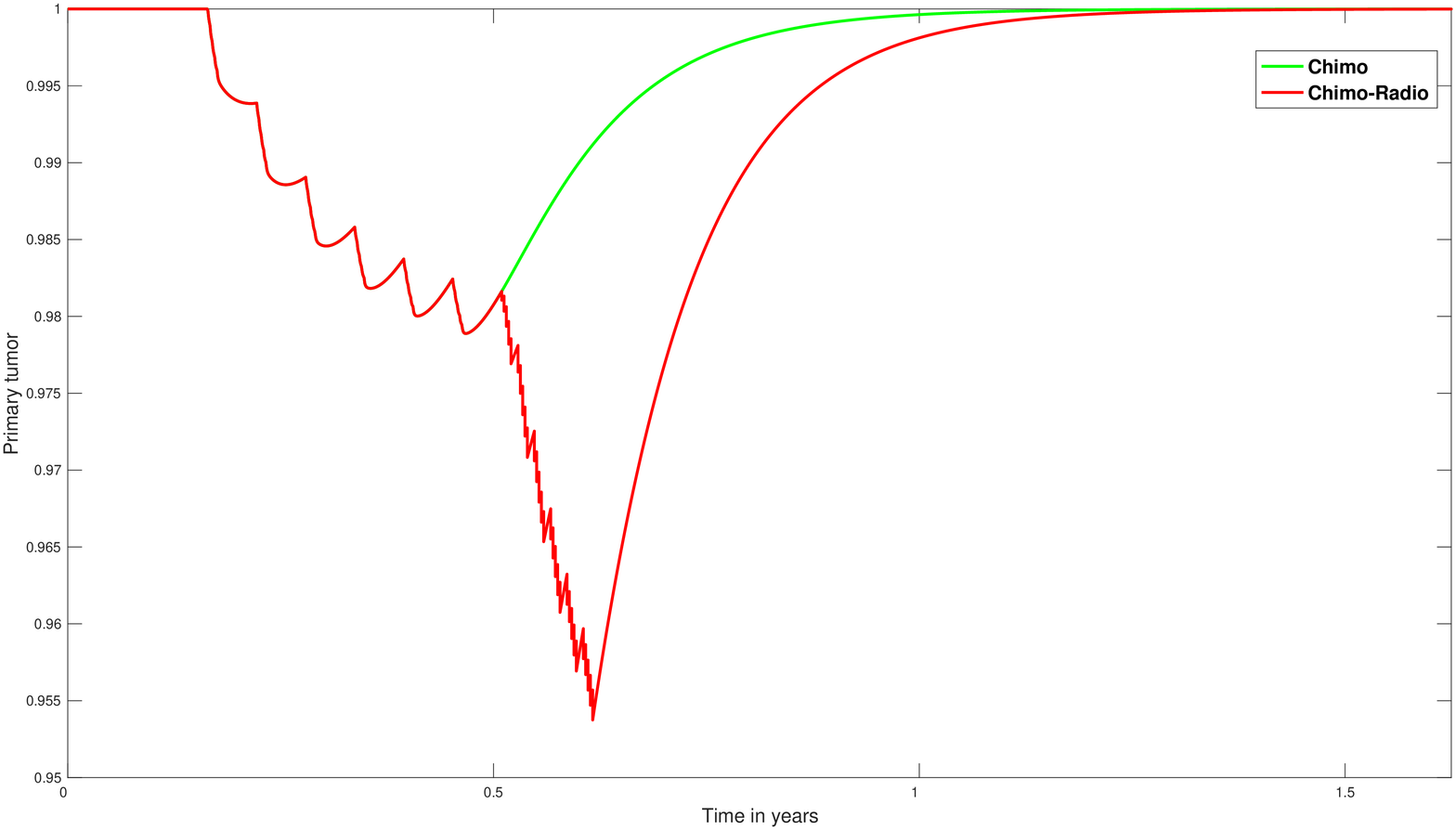}
    \caption{The ratio of the number of cells,}
    \label{fig10}
\end{subfigure}
\caption{Comparison in terms of cellular evolution of the primary tumor between different types of treatment (chemotherapy, chemo-radiotherapy and no treatment).}
\end{figure}

\begin{figure}[h!]
\begin{subfigure}[b]{.49\linewidth}
    \centering
    \includegraphics[width=0.85\linewidth,height=0.45\linewidth]{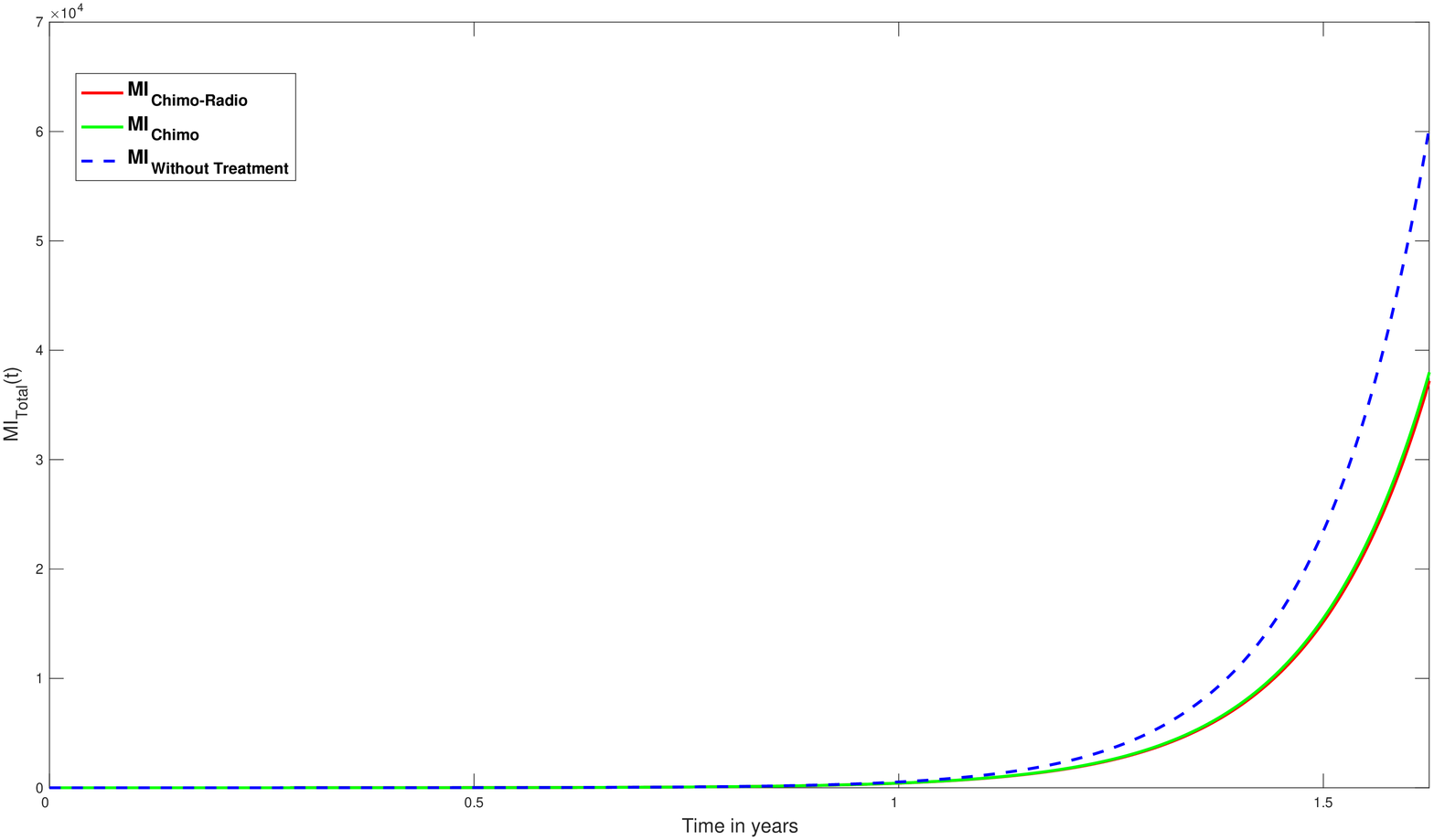}
    \caption{Total MI,}
    \label{fig11}
\end{subfigure}
\hfill
\begin{subfigure}[b]{.45\linewidth}
    \centering
    \includegraphics[width=0.85\linewidth,height=0.45\linewidth]{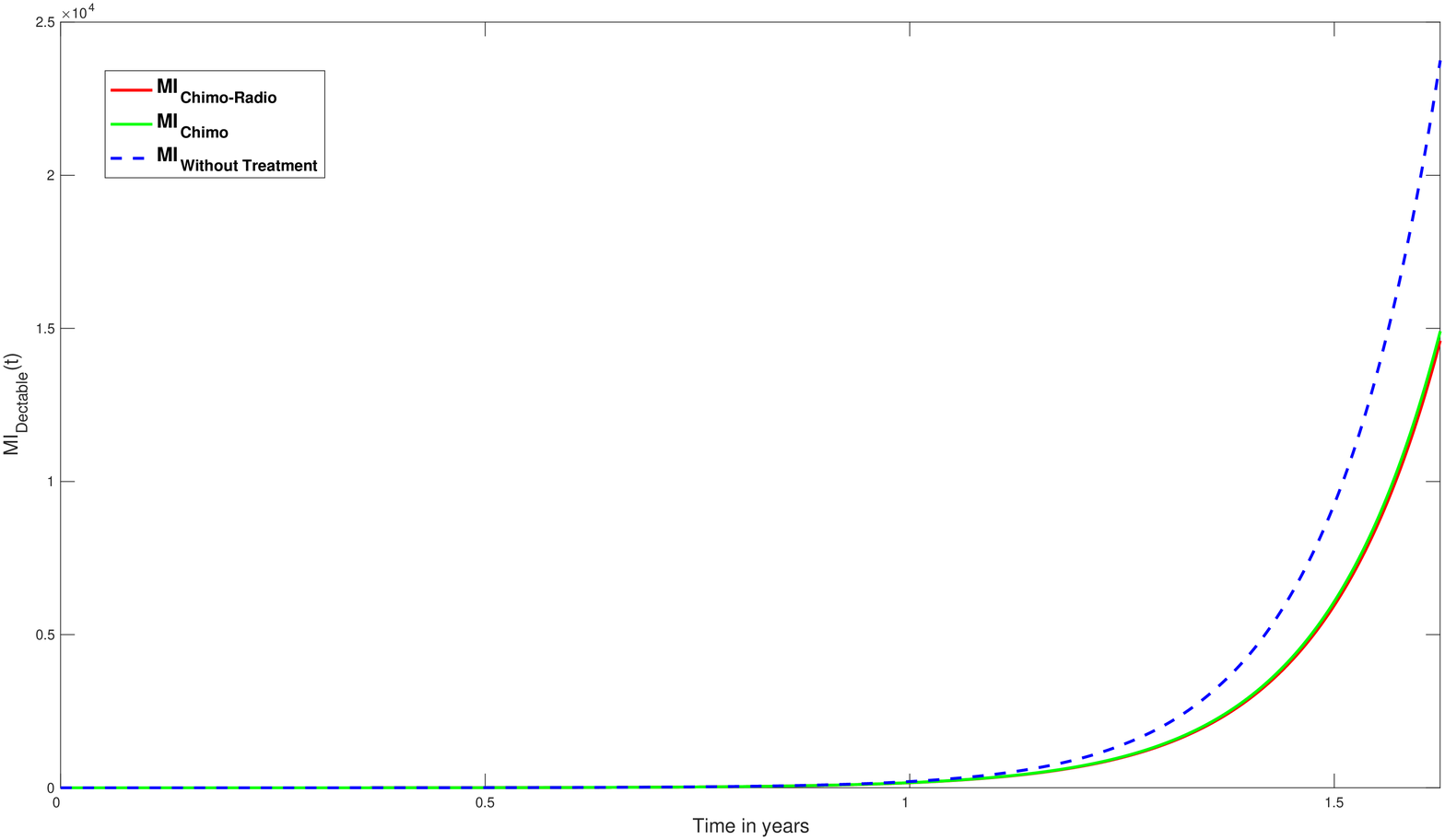}
     \caption{Detectable MI,}
    \label{fig12}
\end{subfigure}
\caption{Comparison in terms of evolution in number of metastases between different types of treatment (chemotherapy, chemo-radiotherapy and no treatment).}
\end{figure}

\begin{figure}[h!]
\begin{subfigure}[b]{.46\linewidth}
    \centering
    \includegraphics[width=0.85\linewidth,height=0.45\linewidth]{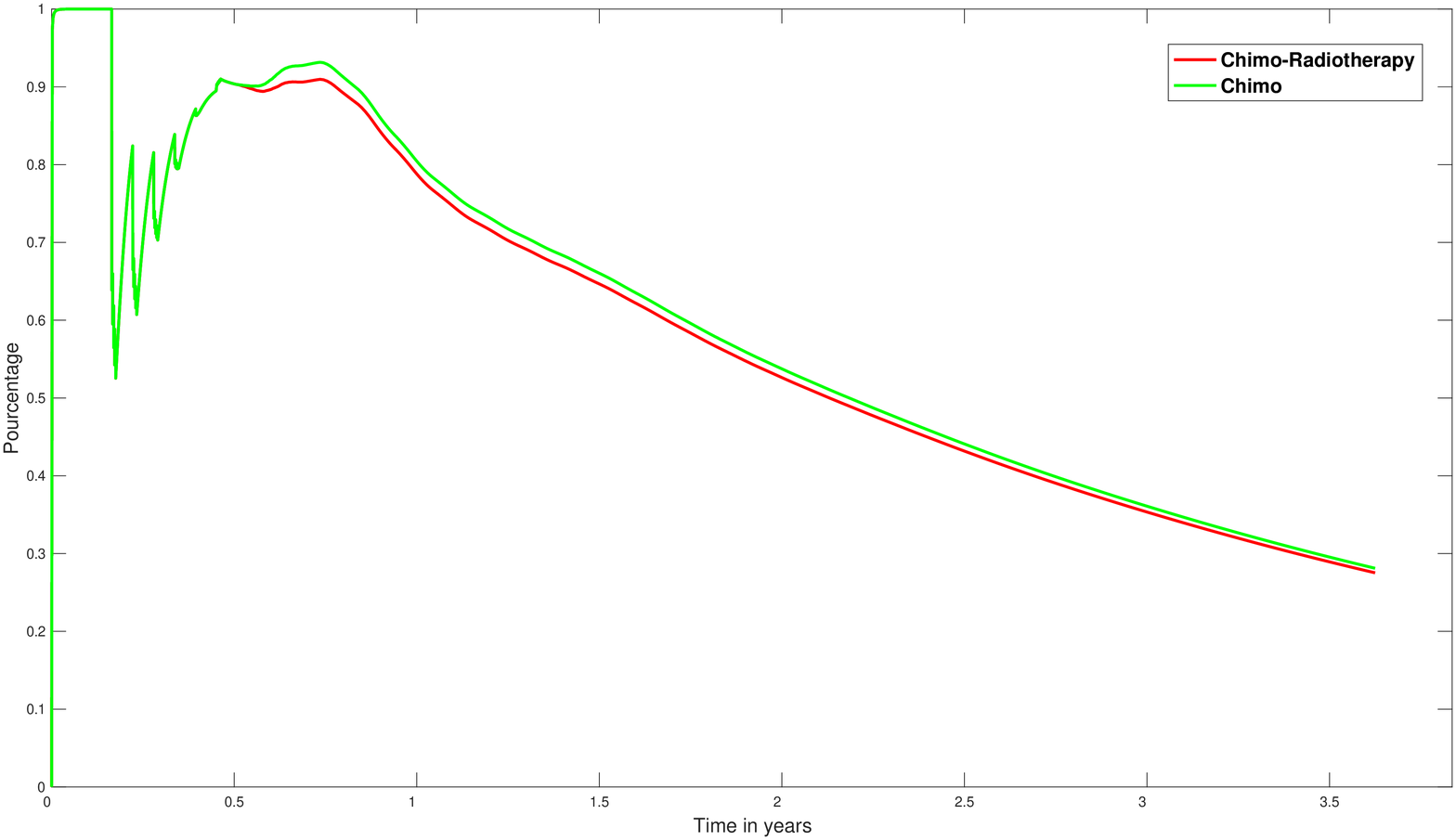}
    \caption{The ratio of the total MI,}
    \label{fig13}
\end{subfigure}
\hfill
\begin{subfigure}[b]{.46\linewidth}
    \centering
    \includegraphics[width=0.85\linewidth,height=0.45\linewidth]{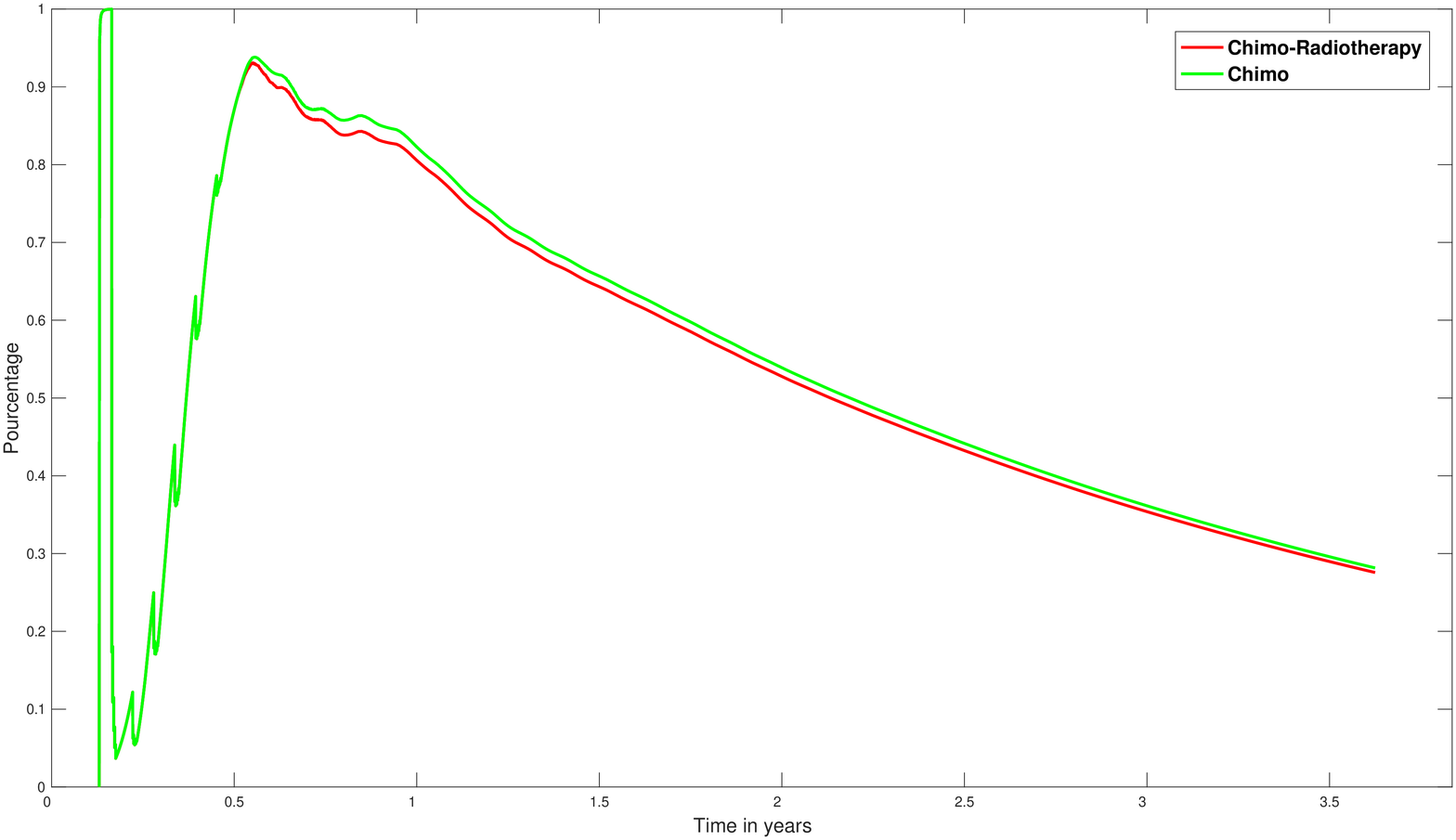}
    \caption{The ratio of the detectable MI,}
    \label{fig14}
\end{subfigure}
\caption{Percent calculation between metastatic index with and without treatment.}
\end{figure}

We consider the effect of concomitant treatments on the evolution, first of a primary tumor, then of a metastatic growing tumor. We consider the following parameters \cite{FV,Faivre}
\begin{equation}\label{parametre}
  a=0.0231,\;\;  b=10^{12},\;\;   m = 5.3\times 10^{- 8},\;\; \alpha=0.55,\;\; \lambda=0.15,\;\; \frac{\lambda}{\mu} = 1.
\end{equation}%
\paragraph{\textbf{Primary tumor}}~\\ 
In Figure \ref{fig9}, we compare the numbers of cells in a primary tumor, with and without treatment. In Figure \ref{fig10}, we present the effect of treatments, in terms of ratio of the cells numbers. The advantage of the mixed treatment is highlighted. 
\paragraph{\textbf{Growing tumor}}~\\
We present, respectively in Figure \ref{fig11} and Figure \ref{fig12}, the effect of treatments on the number of metastatic tumors, detectable and not by medical imaging.
The ratio of metastatic indexes, with and without treatments, is presented in Figure \ref{fig13} and Figure \ref{fig14}. We notice the advantage of the concomitant treatment for our selected  parameters. 

Table \ref{tab2} shows the evolution of the metastatic index with and without treatment at the times of the beginning and the end of treatments and one years later. We can see that one years after the end of treatment, the variation in cell numbers between the chemotherapy alone and the mixed treatment is of the order of $9.9\times 10^4$ cells and $9.8\times 10^4$ cells respectively. We see clearly in Figure \ref{fig13} and Figure \ref{fig14} the effect of concomitant treatment to reduce tumour size.

In conclusion, the metastatic index should enable to better predict the risk of invasive cancer, even it cannot be seen with standard imaging techniques. We notice that the effect of the treatment is strictly linked to the choice of parameters associated to the decease and to the treatments, namely, the pharmacokinetic and pharmacodynamic parameters for chemotherapy and the radiosensitivity parameters for radiotherapy.

\begin{center}
\begin{table}[h!]%
\centering
\begin{tabular}{lccc|}
\cline{2-4}
& \multicolumn{1}{|c}{$t=t_0$}  & $t=T_1$  & \multicolumn{1}{c|}{$t=T$} \\
\hline
\multicolumn{1}{|l}{\textbf{Total Metastatic Index}}&   &  & \\

\multicolumn{1}{|l}{$\text{MI}_{\text{Without Treat}}$}& $0.0413$  & $43.08$   & $1.36\text{e}05$\\
\multicolumn{1}{|l}{$\text{MI}_{\text{Chemo-Radio}}$}  & $0.0440$  & $27.07$   & $3.71\text{e}04$\\
\multicolumn{1}{|l}{$\text{MI}_{\text{Only Chemo}}$}   & $0.0440$  & $27.59$   & $3.79\text{e}04$\\
\hline
\multicolumn{1}{|l}{\textbf{Detectable Metastatic Index }}&   &  &   \\

\multicolumn{1}{|l}{$\text{MI}_{\text{Without Treat}}$}& $3.84\text{e}-6$  & $15.95$   & $4.9\text{e}05$\\
\multicolumn{1}{|l}{$\text{MI}_{\text{Chemo-Radio}}$}  & $3.84\text{e}-6$  & $14.16$   & $3.65\text{e}04$\\
\multicolumn{1}{|l}{$\text{MI}_{\text{Only Chemo}}$}   & $3.84\text{e}-6$  & $14.41$   & $3.74\text{e}04$\\
\hline
\end{tabular}
\caption{Evolution of the metastatic index at the time points: start of treatment ($t_0$), end of treatment ($T_1$) and one years after the end of treatment ($T=T_1+1$).\label{tab2}}
\end{table}
\end{center}

\section{Conclusion}
In this work, We have presented an analysis of a simplified model of tumor growth with treatment. The data associated with the treatment are irregular in time. Our model and existence result could be valid for other problems of population dynamics. We have carried out some numerical tests that could help in the analysis of the treatment effect. Many points remain to be treated. In particular, the inverse problems relating to some parameters of the model are interesting both theoretically, numerically and clinically.

\subsection{Acknowledgement} The authors would like to thank Belhassen Dehman, Professor at the Faculty of Sciences of Tunis, University Tunis El Manar, for the fruitful discussions on the theoretical part of this work.

\printbibliography 

\end{document}